\documentclass[preprint,12pt]{elsarticle}
\usepackage{graphicx}
\usepackage[figuresright]{rotating}
\usepackage{amssymb}
\usepackage{amsthm}
\usepackage{multirow}
\usepackage{multicol}
\usepackage{subfigure}
\usepackage{epstopdf}
\usepackage{fullpage}
\usepackage{latexsym,amsmath}
\usepackage{stmaryrd,appendix}
\usepackage{algorithm,algorithmic}
\usepackage{subfigure}
\usepackage{caption}
\usepackage{url}
\usepackage{amsfonts}
\usepackage{xcolor}
\usepackage{diagbox}
\usepackage{array}
\usepackage{booktabs}
\usepackage{rotating} 
\usepackage{wrapfig}
\biboptions{sort&compress}

\newtheorem{lemma}{Lemma}
\newtheorem{theorem}{Theorem}
\newtheorem{definition}{Definition}
\usepackage{setspace} 
\usepackage{amssymb}
\begin{document}
\begin{spacing}{1.95}
\begin{frontmatter}



\title{Low-rank quaternion tensor completion for recovering color videos and images}

\author[lab1]{Jifei Miao}
\ead{jifmiao@163.com}
\author[lab1]{Kit Ian Kou\corref{cor1}}
\ead{kikou@umac.mo}
\author[lab2]{Wankai Liu}
\ead{zjnulwk@163.com}

\address[lab1]{Department of Mathematics, Faculty of
	Science and Technology, University of Macau, Macau 999078, China}
\address[lab2]{School of Mathematics and Quantitative Economics, Shandong University of Finance and
	Economics, Shandong, China}

\cortext[cor1]{Corresponding author}

\begin{abstract}
	Low-rank quaternion tensor completion method, a novel approach to recovery color videos and images is proposed in this paper. We respectively reconstruct a color image and a color video as a quaternion matrix (second-order tensor) and a third-order quaternion tensor by encoding the red, green, and blue channel pixel values on the three imaginary parts of a quaternion. Different from some traditional models which treat color pixel as a scalar and represent color channels separately, whereas, during the quaternion-based reconstruction, it is significant that the inherent color structures of color images and color videos can be completely preserved. Under the definition of Tucker rank, the global low-rank prior to quaternion tensor is encoded as the nuclear norm of unfolding quaternion matrices. Then, by applying the ADMM framework, we provide the tensor completion algorithm for any order $(\geq 2)$ quaternion tensors, which theoretically can be well used to recover missing entries of any multidimensional data with color structures. Simulation results for color videos and color images recovery show the superior performance and efficiency of the proposed method over some state-of-the-art existing ones.
\end{abstract}

\begin{keyword}
	Quaternion \sep color videos  \sep color images \sep tensor completion \sep low-rank.
\end{keyword}
\end{frontmatter}
%---------------------------------------------------------------------%

\section{Introduction}
\label{section_1}
A color video or color image contains red, blue, and green channels. In most cases, some data of the acquired videos or images are missed during acquisition and transmission, which poses great challenges to further processing of them. Hence, a well-performed recovery technology is important and necessary to recover complete videos or images from their incomplete observations. The core of the missing value estimation lies on how to exactly build a proper low-rank regularizer to measure the global structure of the underlying video or image data according to the fact that the data inherently possess a low-rank structure \cite{DBLP:journals/pr/LiWLT19}.

In the past few decades, the low-rank matrix completion problem has been widely studied and proven very useful
in the application of images and even videos recovery \cite{DBLP:journals/pr/FanC17a, DBLP:journals/pr/FanC18}. Commonly, the method is to stack all the image or video pixels as column vectors of a matrix, and recovery theories and algorithms
are adopted to the resulting matrix which is low-rank or approximately low-rank. However, these image and video recovery models are usually developed for grey-level pixels. For color videos and color images processing, traditional matrix-based methods usually ignore the mutual connection among channels, because these recovery methods are applied to red, green, and blue channels separately, which is likely to result in color distortion during the recovery process \cite{DBLP:journals/tip/ZouKW16}.

On the other hand, with the success of low-rank matrix completion, low-rank tensor completion is an extension to process the multidimensional data \cite{DBLP:journals/pr/QinJHLLZDZF19,DBLP:journals/tip/ZhouLLZ18,DBLP:journals/sigpro/LongLCZ19}. A color image with red, blue, and green channels can be naturally regarded as a third-order tensor. Each frontal slice of this third-order tensor corresponds to a channel of the color image.  Analogously, a video comprised of color images is a fourth-order tensor with an additional index for a temporal variable \cite{DBLP:journals/sigpro/LongLCZ19}. Nevertheless, there are still some underlying restrictions on tensor-based completion algorithms, especially for color video recovery problems. For example there are plenty of  completion algorithms using  tensor singular value decomposition (t-SVD) and the tubal rank \cite{6909886TNN,DBLP:journals/tip/ZhouLLZ18}, \emph{etc.}, can not be well applied to color video recovery problem, since the t-SVD and the tubal rank theories \cite{DBLP:journals/siammax/KilmerBHH13} they are based on are defined for third-order tensors. Therefore, for grey scalar videos (third-order tensors), they can obtain well performance, however, for color videos (fourth-order tensors),  they may ignore the inherent color structures, and then can not offer a satisfying recovery result. In addition, there are some factorization based approaches for tensor completion, such as CANDECOMP/PARAFAC (CP) and Tucker factorizations based approaches \cite{DBLP:journals/pami/ZhaoZC15,DBLP:journals/tsp/YokotaZC16}, tensor unfolding based approaches \cite{DBLP:journals/pami/LiuMWY13,1930-8337_2015_2_601}, \emph{etc.}. These
methods can deal with color videos directly, however, the factorization or matricization operation may destroy color pixel structure and lead to color distortion during the recovery process \cite{DBLP:journals/pr/YuWGGXP19}.

\begin{figure}[htbp]
\centering
\includegraphics[width=3cm,height=2.2cm]{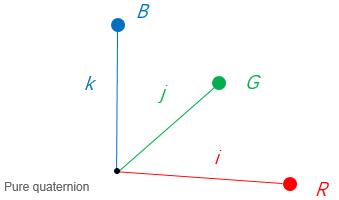}
\caption{The graphical of a pure quaternion representing a color pixel. $i$, $j$ and $k$ correspond to three channels ($R$, $G$ and $B$) of the color pixel.}
\label{quaternion_fig}
\end{figure}

Different from conventional matrix and tensor based models, in this paper, we make use of quaternion tensors to represent color videos and color images\footnote{Color images are represented by second-order quaternion tensors, we also call them quaternion matrices in the paper.}, and study the problem of quaternion tensors completion to estimate missing data of them. Actually, the quaternion has achieved excellent results in color image processing problems including histopathological image analysis \cite{DBLP:journals/pr/ShiZWGZY19}, color image denoising and representation \cite{DBLP:journals/mssp/GaiYW015, DBLP:journals/pr/HosnyD19}, color object recognition \cite{DBLP:journals/pr/LiLHS15, DBLP:journals/pr/ShaoSWCC14}, and so on. The quaternion based method is to encode the red, green, and blue channel pixel values on the three imaginary parts of a quaternion \cite{DBLP:journals/pr/LiLHS15}. That is
\begin{equation}
\label{equ1}
\dot{t}=0+t_{r}i+t_{g}j+t_{b}k,
\end{equation}
where $\dot{t}$ denotes a color pixel, $t_{r}$, $t_{g}$ and $t_{b}$ are, respectively, the red, green and blue channel pixel values, and $i$, $j$ and $k$ are the three imaginary units\footnote{A detailed introduction to the quaternion can be found in Section \ref{section_2}.}. The graphical of a pure quaternion representing a color pixel can be seen in Figure.\ref{quaternion_fig}.
By using (\ref{equ1}), an $M\times N\times3$ color
image, and an $M\times N\times3\times T$ (image row $\times$ image column $\times$ RGB $\times$ frame) color video are respectively described by a second-order quaternion tensor (quaternion matrix) with size $M\times N$ and a third-order quaternion tensor with size $M\times N \times T$ whose entries are pure quaternions. The main advantage of this quaternion representation is that it processes a color pixel holistically as a vector field and handles the coupling between the color channels naturally \cite{ DBLP:journals/pr/HosnyD19}, and color information of source video and image is fully utilized.

However, the existing quaternion based methods mainly consider the color image issues but not consider the higher dimensional data structures, for example, color videos.  In this paper, we reconstruct an $M\times N\times3\times T$ color video (fourth-order tensor) as an $M\times N \times T$ third-order quaternion tensor with each frontal slice being a quaternion matrix. It is important to highlight that different from traditional tensor model, the three color channels of each pixel in quaternion tensor  can be fully connected by the model (\ref{equ1}), and it is clear that even if the matricization operations (\emph{e.g.}, tensor unfolding) can not destroy the color pixel structure, \emph{i.e.,} the relative positions of the three color channel pixels of one pixel will remain unchanged just as Figure.\ref{quaternion_fig}. Hence, Under the definition of Tucker rank, the global low-rank prior to quaternion tensor is encoded as the nuclear norm of unfolding quaternion matrices. Then, we provide the completion algorithm for any order $(\geq 2)$ quaternion tensors by applying the alternating direction method of multipliers (ADMM) \cite{DBLP:journals/ftml/BoydPCPE11} framework. Simulation results for color videos and color images recovery show that
the performance of the proposed method is better than that of the testing methods.

The rest of this paper is organized as follows. Section \ref{section_2} introduces some notations and preliminaries for quaternion algebra and quaternion tensor. Section \ref{sec3} reviews the tensor completion theory and proposes our quaternion-based tensor completion model. The detailed overview of the quaternion tensor completion method is presented in Section \ref{sec4}. Section \ref{sec5} provides some simulations to illustrate the performance of our approach, and compare it with some state-of-the-art methods. Finally, some conclusions are drawn in Section \ref{sec6}.

\section{Notations and preliminaries}
\label{section_2}
In this section, we first summarize some main notations and then introduce some basic knowledge of quaternion algebra and quaternion tensor.

\subsection{Notations and definitions}
In this paper, $\mathbb{R}$, $\mathbb{C}$, and $\mathbb{H}$ respectively denote the set of real numbers, the set of complex numbers and the set of quaternions. A scalar, a vector, a matrix, and a tensor are written as $a$, $\mathbf{a}$, $\mathbf{A}$, and $\mathcal{A}$, respectively. For
a tensor $\mathcal{A}$, we use the Matlab notation $\mathcal{A}(:,:,k)$  to denote its $k$-th frontal slice and the $\mathbf{A}_{[k]}$ to denote its mode-k unfolding. A dot (above the variable) is used to denote a quaternion variable (\emph{e.g.,} \cite{DBLP:journals/tip/ZouKW16, DBLP:journals/tip/XuYXZN15}), $\dot{a}$, $\dot{\mathbf{a}}$, $\dot{\mathbf{A}}$ and $\dot{\mathcal{A}}$
respectively represent a quaternion scalar, a quaternion vector, a quaternion matrix and a quaternion tensor. $(\cdot)^{\ast}$ and $(\cdot)^{H}$ denote the
conjugation and conjugate transpose.  $|\cdot|$, $\|\cdot\|_{F}$ and $\|\cdot\|_{\ast}$ are respectively the moduli, the Frobenius norm, and the nuclear norm. ${\rm{tr}}\{\cdot\}$, ${\rm{rank}}(\cdot)$ and $\nabla_{sub}$ denote the trace, rank and subgradient operators respectively. ${\rm{Unfold}}_{k}$ is the Mode-k unfolding operator of tensors, and we use ${\rm{Fold}}_{k}$ to denote the inverse operator of ${\rm{Unfold}}_{k}$.

\subsection{Basic knowledge of quaternion algebra}
Quaternions were discovered in 1843 by W.R. Hamilton \cite{doi:10.1080/14786444408644923}\footnote{Here we just give some fundamental algebraic operations used in our work briefly, which follow the definition in \cite{DBLP:journals/pr/ShaoSWCC14,Girard2007Quaternions}.  Readers
can find more details on quaternion algebra in the references.}. A quaternion $\dot{q}\in\mathbb{H}$  is a four-dimensional (4D) hypercomplex number and has a Cartesian form given by:
\begin{equation}
\label{equ2}
\dot{q}=q_{0}+q_{1}i+q_{2}j+q_{3}k,
\end{equation}
where $q_{l}\in\mathbb{R}\: (l=0,1,2,3)$ are called its components, and $i, j, k$ are
square roots of -1 and are related through the famous relations:
\begin{align}
\left\{
\begin{array}{lc}
i^{2}=j^{2}=k^{2}=ijk=-1,\\
ij=-ji=k,
jk=-kj=i, 
ki=-ik=j.
\end{array}
\right.
\end{align}
A quaternion $\dot{q}\in\mathbb{H}$ can be decomposed into a real part R$(\dot{q})$ and an imaginary part I$(\dot{q})$:
\begin{equation}
\label{equ3}
\dot{q}={\rm{R}}(\dot{q})+{\rm{I}}(\dot{q}),
\end{equation}
where ${\rm{R}}(\dot{q})=q_{0}$, ${\rm{I}}(\dot{q})=q_{1}i+q_{2}j+q_{3}k$. Then, $\dot{q}\in\mathbb{H}$ will be called a
pure quaternion if its real part is null, \emph{i.e.}, if ${\rm{R}}(\dot{q})=0$.
Given two quaternions $\dot{p}$ and $\dot{q}\in\mathbb{H}$, the sum and multiplication of them are respectively:
\begin{equation*}
\label{equ4}
\dot{p}+\dot{q}=(p_{0}+q_{0})+(p_{1}+q_{1})i+(p_{2}+q_{2})j+(p_{3}+q_{3})k
\end{equation*}
and
\begin{align*}
\label{equ5}
\dot{p}\dot{q}=&(p_{0}q_{0}-p_{1}q_{1}-p_{2}q_{2}-p_{3}q_{3})
+(p_{0}q_{1}+p_{1}q_{0}+p_{2}q_{3}-p_{3}q_{2})i\\
&+(p_{0}q_{2}-p_{1}q_{3}+p_{2}q_{0}+p_{3}q_{1})j
+(p_{0}q_{3}+p_{1}q_{2}-p_{2}q_{1}+p_{3}q_{0})k.
\end{align*}
It is noticeable that the multiplication of two quaternions is not
commutative so that in general $\dot{p}\dot{q}\neq\dot{q}\dot{p}$.

The conjugate and the modulus of a quaternion $\dot{q}\in\mathbb{H}$ are,
respectively, defined as follows:
\begin{align*}
\dot{q}^{\ast}=q_{0}-q_{1}i-q_{2}j-q_{3}k,\\
|\dot{q}|=\sqrt{q_{0}^{2}+q_{1}^{2}+q_{2}^{2}+q_{3}^{2}}.
\end{align*}
Every quaternion $\dot{q}\in\mathbb{H}$ can be uniquely represented as the Cayley–Dickson (CD) form:
\begin{equation}
\label{equ6}
\dot{q}=z_{1}+z_{2}j,
\end{equation}
where $z_{1}=q_{0}+q_{1}i$ and $z_{2}=q_{2}+q_{3}i$ are complex numbers.

\subsection{Quaternion matrix and tensor}
The quaternion matrix is denoted as $\dot{\mathbf{Q}}=(\dot{q}_{n_{1},n_{2}})\in\mathbb{H}^{N_{1}\times N_{2}}$, $1\leq n_{1}\leq N_{1}$, $1\leq n_{2}\leq N_{2}$, where each entry is a
quaternion \cite{ZHANG199721}. We often rewritten it as
\begin{equation*}
\dot{\mathbf{Q}}=\mathbf{Q}_{0}+\mathbf{Q}_{1}i+\mathbf{Q}_{2}j+\mathbf{Q}_{3}k,
\end{equation*}
where $\mathbf{Q}_{l}\in\mathbb{R}^{N_{1}\times N_{2}}\: (l=0,1,2,3)$, $\dot{\mathbf{Q}}$ is named a pure quaternion matrix when $\mathbf{Q}_{0}=\mathbf{0}$. In addition to scalar representations, based on the CD form (\ref{equ6}), there exists their isomorphic complex representation denoted as $f(\dot{\mathbf{Q}})\in\mathbb{C}^{2N_{1}\times 2N_{2}}$,  is of the form:
\begin{equation}
\label{definef}
f(\dot{\mathbf{Q}})=\left(\begin{array}{cc}
\mathbf{Z}_{1}	&  \mathbf{Z}_{2}\\
-\mathbf{Z}^{\ast}_{2}	& \mathbf{Z}^{\ast}_{1}
\end{array}\right),
\end{equation}
where $\mathbf{Z}_{1}=\mathbf{Q}_{0}+\mathbf{Q}_{1}i$, $\mathbf{Z}_{2}=\mathbf{Q}_{2}+\mathbf{Q}_{3}i\in\mathbb{C}^{N_{1}\times N_{2}}$.

Note that the multiplication between quaternion matrices can be defined similar to classical multiplication between real or complex matrices, except that the multiplication between two quaternion numbers is employed.
\begin{definition}(The rank of quaternion matrix \cite{ZHANG199721})
The maximum number of right (left) linearly independent columns (rows) of a quaternion matrix $\dot{\mathbf{Q}}\in\mathbb{H}^{N_{1}\times N_{2}}$ is called the
rank of $\dot{\mathbf{Q}}$.
\end{definition}

\begin{theorem}(Quaternion singular value decomposition (QSVD) \cite{ZHANG199721})
Let $\dot{\mathbf{Q}}\in\mathbb{H}^{N_{1}\times N_{2}}$ be of rank r.  Then, there exist two unitary quaternion matrices\footnote{A unitary quaternion matrix $\dot{\mathbf{U}}\in\mathbb{H}^{N\times N}$ has the following
	property: $\dot{\mathbf{U}}\dot{\mathbf{U}}^{H}=\dot{\mathbf{U}}^{H}\dot{\mathbf{U}}=\mathbf{I}_{N}$, with $\mathbf{I}_{N}\in\mathbb{R}^{N\times N}$ being the quaternion identity matrix which is the same as the
	classical identity matrix.} $\dot{\mathbf{U}}\in\mathbb{H}^{N_{1}\times N_{1}}$ and $\dot{\mathbf{V}}\in\mathbb{H}^{N_{2}\times N_{2}}$ such that	
\begin{equation*}
\dot{\mathbf{U}}^{H}\dot{\mathbf{Q}}\dot{\mathbf{V}}=\mathbf{\Lambda}=\left(\begin{array}{cc}
\mathbf{\Sigma}_{r}	&  \mathbf{0}\\
\mathbf{0}	& \mathbf{0}
\end{array}\right),
\end{equation*}
where $\mathbf{\Sigma}_{r}={\rm{diag}}(\sigma_{1},\ldots,\sigma_{r})$ is a real diagonal matrix and has $r$  positive entries $\sigma_{k}, \,(k=1,\ldots,r)$ on its diagonal (\emph{i.e.}, positive singular values of $\dot{\mathbf{Q}}$).
\end{theorem}
The relation between the QSVD of quaternion matrix  $\dot{\mathbf{Q}}\in\mathbb{H}^{N_{1}\times N_{2}}$ and the SVD of its isomorphic complex matrix $f(\dot{\mathbf{Q}})\in\mathbb{C}^{2N_{1}\times 2N_{2}}$ ($f(\dot{\mathbf{Q}})=\mathbf{U}\check{\mathbf{\Lambda}}\mathbf{V}^{H}$) is defined as \cite{DBLP:journals/tip/XuYXZN15}:
\begin{align}
\label{qsvd}
\left\{
\begin{array}{lc}
\mathbf{\Lambda}={\rm{row}}_{odd}({\rm{col}}_{odd}(\check{\mathbf{\Lambda}})),\\
\dot{\mathbf{U}}={\rm{col}}_{odd}(\mathbf{U}_{1})+{\rm{col}}_{odd}(-(\mathbf{U}_{2})^{\ast})j,\\
\dot{\mathbf{V}}={\rm{col}}_{odd}(\mathbf{V}_{1})+{\rm{col}}_{odd}(-(\mathbf{V}_{2})^{\ast})j,
\end{array}
\right.
\end{align}
such that $\dot{\mathbf{Q}}=\dot{\mathbf{U}}\mathbf{\Lambda}\dot{\mathbf{V}}^{H}$, where
\begin{align*}
\mathbf{U}=\left(\begin{array}{c}
(\mathbf{U}_{1})_{N_{1}\times 2N_{1}} \\
(\mathbf{U}_{2})_{N_{1}\times 2N_{1}}
\end{array} \right),\quad
\mathbf{V}=\left(\begin{array}{c}
(\mathbf{V}_{1})_{N_{2}\times 2N_{2}} 	\\
(\mathbf{V}_{2})_{N_{2}\times 2N_{2}}
\end{array} \right),
\end{align*}
and ${\rm{row}}_{odd}(\mathbf{M})$, ${\rm{col}}_{odd}(\mathbf{M})$ respectively extract the odd rows and odd columns of matrix $\mathbf{M}$. Based on the QSVD, we define the quaternion matrix  nuclear norm (QMNN) below.
\begin{definition}(QMNN) Given a quaternion matrix $\dot{\mathbf{Q}}\in\mathbb{H}^{N_{1}\times N_{2}}$, the QMNN of it is defined as
\begin{equation}
\label{qmnn}
\|\dot{\mathbf{Q}}\|_{\ast}=\sum_{k=1}^{{\rm{min}}\{N_{1},N_{2}\}}\sigma_{k},
\end{equation}	
where $\sigma_{k}$ is the singular value of $\dot{\mathbf{Q}}$, which can be obtained by the QSVD of $\dot{\mathbf{Q}}$.
\end{definition}
In addition, the Frobenius norm of the quaternion matrix $\dot{\mathbf{Q}}\in\mathbb{H}^{N_{1}\times N_{2}}$ is defined as \cite{ZHANG199721}:
$\|\dot{\mathbf{Q}}\|_{F}=\sqrt{\sum_{n_{1}=1}^{N_{1}}\sum_{n_{2}=1}^{N_{2}}|\dot{q}_{n_{1},n_{2}}|^{2}}=\sqrt{{\rm{tr}}\{(\dot{\mathbf{Q}})^{H}\dot{\mathbf{Q}}\}}$.

Analogously, in this paper, we generalize the definition of quaternion matrix to higher dimensional quaternion array, \emph{i.e.}, quaternion tensor.
\begin{definition}(Quaternion tensor) A multidimensional array or an Nth-order tensor is called a quaternion tensor if its entries are quaternion numbers, \emph{i.e.},
\begin{align}
\label{qtensor1}
\dot{\mathcal{T}}&=(\dot{t}_{n_{1},n_{2},\ldots,n_{N}})\in\mathbb{H}^{N_{1}\times N_{2} \times\ldots \times N_{N}}	\nonumber\\
&=\mathcal{T}_{0}+\mathcal{T}_{1}i+\mathcal{T}_{2}j+\mathcal{T}_{3}k,
\end{align}	
where $\mathcal{T}_{l}\in\mathbb{R}^{N_{1}\times N_{2} \times\ldots \times N_{N}}\: (l=0,1,2,3)$, $\dot{\mathcal{T}}$ is named a pure quaternion tensor when $\mathcal{T}_{0}$ is a zero tensor.
\end{definition}

\begin{definition}(Mode-k unfolding) For an Nth-order quaternion tensor $\dot{\mathcal{T}}\in\mathbb{H}^{N_{1}\times N_{2} \times\ldots \times N_{N}}$, its mode-k unfolding is defined as a quaternion matrix
\begin{equation*}
{\rm{Unfold}}_{k}(\dot{\mathcal{T}})=\dot{\mathbf{T}}_{[k]}\in\mathbb{H}^{N_{k}\times N_{1}\ldots N_{k-1}N_{k+1}\ldots N_{N}}\quad  \text{with entries}
\end{equation*}
\begin{equation*}
\dot{\mathbf{T}}_{[k]}(n_{k}, n_{1}\ldots n_{k-1}n_{k+1}\ldots n_{N})=\dot{t}_{n_{1},n_{2},\ldots,n_{N}},
\end{equation*}
where $\dot{t}_{n_{1},n_{2},\ldots,n_{N}}$ is the $(n_{1},n_{2},\ldots,n_{N})$th entry of $\dot{\mathcal{T}}$.
\end{definition}

\begin{definition}(Tucker rank \cite{DBLP:journals/siamrev/KoldaB09}) Given a quaternion tensor $\dot{\mathcal{T}}\in\mathbb{H}^{N_{1}\times N_{2} \times\ldots \times N_{N}}$, the Tucker rank of it is defined as
\begin{equation}
\label{Tucker_rank}
{\rm{rank_{tucker}}}(\dot{\mathcal{T}})=({\rm{rank}}(\dot{\mathbf{T}}_{[1]}),{\rm{rank}}(\dot{\mathbf{T}}_{[2]}),\ldots,{\rm{rank}}(\dot{\mathbf{T}}_{[N]})),	
\end{equation}
where ${\rm{rank}}(\dot{\mathbf{T}}_{[k]})$ denotes the rank of the mode-k unfolding quaternion matrix $\dot{\mathbf{T}}_{[k]}$.
\end{definition}

\section{Problem formulation}
\label{sec3}
In this section, we first review the tensor completion theory and then propose our quaternion based tensor
completion model.

\subsection{Tensor completion theory}
The tensor completion problem consists of recovering a tensor from a subset of its entries. The key is to build up the relationship between the available and the missing entries \cite{DBLP:journals/pami/LiuMWY13}. The usual structural assumption on a tensor that makes the problem well-posed is that the tensor is low-rank or approximate low-rank. Mathematically, the optimization model for tensor completion problem can be formulated as:
\begin{equation}
\label{equ7}
\begin{split}
&\mathop{{\rm{minimize}}}\limits_{\mathcal{T}}\quad {\rm{rank}}(\mathcal{T})\\
&\text{subject to} \quad  \mathcal{P}_{\Omega}(\mathcal{T})=\mathcal{Y},
\end{split}
\end{equation}
where $\mathcal{Y}$ is the underlying complete tensor, $\mathcal{T}$  is the observed tensor, and $\mathcal{P}_{\Omega}$ denotes the random sampling operator which is defined by:

\begin{equation*}
\mathcal{P}_{\Omega}(\mathcal{T})=\left\{
\begin{array}{lc}
 t_{n_{1},n_{2},\ldots,n_{N}},\qquad &(n_{1},n_{2},\ldots,n_{N})\in \Omega, \\
0,  &\text{otherwise}.
\end{array}
\right.
\end{equation*}
However, there is no unique definition for the rank of tensors, such as CP rank \cite{doi:10.1137/0614071}, Tucker rank \cite{DBLP:journals/pami/LiuMWY13}, tubal rank \cite{DBLP:journals/tip/ZhouLLZ18}, tensor train rank \cite{DBLP:journals/tip/BenguaPTD17}, \emph{etc.}.
With different definitions of tensor rank, there are many methods optimization models for tensor completion problem. Among all definitions for the rank of tensors, the Tucker rank is widely used to depict the low-rankness of the underlying tensor. based on minimizing Tucker rank, (\ref{equ7}) can be formulated as:
\begin{equation}
\label{equ8}
\begin{split}
&\mathop{{\rm{minimize}}}\limits_{\mathcal{T}}\quad {\rm{rank_{tucker}}}(\mathcal{T})\\
&\text{subject to} \quad  \mathcal{P}_{\Omega}(\mathcal{T})=\mathcal{Y},
\end{split}
\end{equation}
According to the definition of Tucker rank, (\ref{equ8}) can be written as \cite{DBLP:journals/pami/LiuMWY13,DBLP:journals/ijon/TanCWZR14}:
\begin{equation}
\label{equ9}
\begin{split}
&\mathop{{\rm{minimize}}}\limits_{\mathbf{T}_{[n]}}\quad \sum_{n=1}^{N}\alpha_{n}{\rm{rank}}(\mathbf{T}_{[n]})\\
&\text{subject to} \quad  \mathcal{P}_{\Omega}(\mathcal{T})=\mathcal{Y},
\end{split}
\end{equation}
where $\alpha_{n}$ are nonnegative constants. However, directly optimizing the problem (\ref{equ9}) is NP-hard \cite{DBLP:journals/siammax/GillisG11}. Inspired by matrix nuclear norm, the tightest convex surrogate of the matrix rank, Liu et al. \cite{DBLP:journals/pami/LiuMWY13} established the following definition of the nuclear norm
for tensors:
\begin{equation}
\label{equ10}
\|\mathcal{T}\|_{\ast}=\sum_{n=1}^{N}\alpha_{n}\|\mathbf{T}_{[n]}\|_{\ast}.
\end{equation}
Then, based on (\ref{equ10}), the problem (\ref{equ9}) can be finally rewritten as:
\begin{equation}
\label{equ11}
\begin{split}
&\mathop{{\rm{minimize}}}\limits_{\mathbf{T}_{[n]}}\quad \sum_{n=1}^{N}\alpha_{n}\|\mathbf{T}_{[n]}\|_{\ast}\\
&\text{subject to} \quad  \mathcal{P}_{\Omega}(\mathcal{T})=\mathcal{Y}.
\end{split}
\end{equation}

\subsection{Proposed formulation of quaternion tensor completion}
Quaternion tensor completion can be regarded as the generalization of the traditional tensor completion problem in the quaternion number field, which is to estimate the missing values of a quaternion tensor $\dot{\mathcal{T}}\in\mathbb{H}^{N_{1}\times N_{2} \times\ldots \times N_{N}}$
under a given subset $\Omega$ of its entries $\{\dot{\mathcal{T}}_{n_{1},n_{2},\ldots,n_{N}}|(n_{1},n_{2},\ldots,n_{N})\in\Omega \}$. That is
\begin{equation}
\label{equ12}
\begin{split}
&\mathop{{\rm{minimize}}}\limits_{\dot{\mathcal{T}}}\quad {\rm{rank}}(\dot{\mathcal{T}})\\
&\text{subject to} \quad  \mathcal{P}_{\Omega}(\dot{\mathcal{T}})=\dot{\mathcal{Y}},
\end{split}
\end{equation}
where $\dot{\mathcal{Y}}$ is the underlying complete quaternion tensor, $\dot{\mathcal{T}}$ is the observed quaternion tensor.

Based on the previous \textbf{Definition 3}, \textbf{Definition 4} and \textbf{Definition 5}, and followed by traditional tensor case, we finally translate the problem (\ref{equ12}) into the following low-rank quaternion tensor completion formulation:
\begin{equation}
\label{equ13}
\begin{split}
&\mathop{{\rm{minimize}}}\limits_{\dot{\mathbf{T}}_{[n]}}\quad \sum_{n=1}^{N}\alpha_{n}\|\dot{\mathbf{T}}_{[n]}\|_{\ast}\\
&\text{subject to} \quad  \mathcal{P}_{\Omega}(\dot{\mathcal{T}})=\dot{\mathcal{Y}}.
\end{split}
\end{equation}
The formulation (\ref{equ13}) can be well used to recover missing entries of any multidimensional data with color structures. For special cases, \emph{i.e.}, $N=3$ and $N=2$, we can deal with color videos and color images recovery problems.
It is important to notice that in these kinds of applications the formulation (\ref{equ13}) outperforms that of (\ref{equ11}), because for traditional tensor the Mode-k unfolding operation may destroy the color pixel structure, but for quaternion tensor, the inherent color structures can be completely preserved during this process.

\section{Proposed algorithm}
\label{sec4}
In this section, we show how to solve the optimization
problem (\ref{equ13}), then we provide complexity analyses of the proposed method.
\subsection{Optimization Procedure}
The optimization problem (\ref{equ13}) can be solved by various methods. For efficiency, we adopt the ADMM framework in this paper, which can support the convergence of the algorithm \cite{DBLP:journals/ftml/BoydPCPE11}. By using additional quaternion tensor $\dot{\mathcal{X}}\in\mathbb{H}^{N_{1}\times N_{2} \times\ldots \times N_{N}}$ with $N$ Mode-k unfolding quaternion matrices
$\dot{\mathbf{X}}_{[1]},\dot{\mathbf{X}}_{[2]},\ldots,\dot{\mathbf{X}}_{[N]}$, we first convert (\ref{equ13}) to the following equivalent problem:
\begin{equation}
\label{equ14}
\begin{split}
\mathop{{\rm{minimize}}}\limits_{\dot{\mathcal{T}},\dot{\mathbf{X}}_{[n]}}\quad &\sum_{n=1}^{N}\alpha_{n}\|\dot{\mathbf{X}}_{[n]}\|_{\ast}\\
\text{subject to} \quad &\dot{\mathbf{T}}_{[n]}=\dot{\mathbf{X}}_{[n]}, \ \text{for all} \ n=1,2,\ldots,N\\
&\mathcal{P}_{\Omega}(\dot{\mathcal{T}})=\dot{\mathcal{Y}}.
\end{split}
\end{equation}
This problem can be solved by the ADMM framework, which minimizes the following augmented Lagrangian function:
\begin{align}
\label{equ15}
&\mathcal{L}(\dot{\mathcal{T}}, \{\dot{\mathbf{X}}_{[n]}\}_{n=1}^{N},\{\dot{\mathbf{F}}_{[n]}\}_{n=1}^{N},\{\beta_{n}\}_{n=1}^{N})\nonumber \\
&=\sum_{n=1}^{N}\left(\alpha_{n}\|\dot{\mathbf{X}}_{[n]}\|_{\ast}+\langle\dot{\mathbf{F}}_{[n]},\dot{\mathbf{T}}_{[n]}-\dot{\mathbf{X}}_{[n]}\rangle+\frac{\beta_{n}}{2}\|\dot{\mathbf{T}}_{[n]}-\dot{\mathbf{X}}_{[n]}\|_{F}^{2}\right),
\end{align}
where $\{\beta_{n}\}_{n=1}^{N}$ are the penalty parameters, $\{\dot{\mathbf{F}}_{[n]}\}_{n=1}^{N}$ are the Lagrange multipliers, which are $N$ Mode-k unfolding quaternion matrices of quaternion tensor $\dot{\mathcal{F}}\in\mathbb{H}^{N_{1}\times N_{2} \times\ldots \times N_{N}}$.

It is clear that although the objective function of (\ref{equ15}) is not jointly convex for all variables, it is convex concerning each variable independently. Hence, A natural way to solve the problem is to iteratively optimize
the augmented Lagrangian function (\ref{equ15}) over one variable, while fixing the others. To update each variable, in the $\tau+1$th iteration, perform the following steps:
\begin{itemize}
\item \textbf{Step 1:} $\dot{\mathcal{T}}^{(\tau+1)}=\mathop{{\rm{arg\, min}}}\limits_{\dot{\mathcal{T}}}\:\mathcal{L}\left(\dot{\mathcal{T}}, (\{\dot{\mathbf{X}}_{[n]}\}_{n=1}^{N})^{(\tau)},(\{\dot{\mathbf{F}}_{[n]}\}_{n=1}^{N})^{(\tau)},(\{\beta_{n}\}_{n=1}^{N})^{(\tau)}\right)$,
\item \textbf{Step 2:} $(\{\dot{\mathbf{X}}_{[n]}\}_{n=1}^{N})^{(\tau+1)}=\mathop{{\rm{arg\, min}}}\limits_{\{\dot{\mathbf{X}}_{[n]}\}_{n=1}^{N}}\:\mathcal{L}\left(\dot{\mathcal{T}}^{(\tau+1)}, \{\dot{\mathbf{X}}_{[n]}\}_{n=1}^{N},(\{\dot{\mathbf{F}}_{[n]}\}_{n=1}^{N})^{(\tau)},(\{\beta_{n}\}_{n=1}^{N})^{(\tau)}\right)$,
\item \textbf{Step 3:} $(\{\dot{\mathbf{F}}_{[n]}\}_{n=1}^{N})^{(\tau+1)}=\mathop{{\rm{arg\, min}}}\limits_{\{\dot{\mathbf{F}}_{[n]}\}_{n=1}^{N}}\:\mathcal{L}\left(\dot{\mathcal{T}}^{(\tau+1)}, (\{\dot{\mathbf{X}}_{[n]}\}_{n=1}^{N})^{(\tau+1)},\{\dot{\mathbf{F}}_{[n]}\}_{n=1}^{N},(\{\beta_{n}\}_{n=1}^{N})^{(\tau)}\right)$,
\item \textbf{Step 4:} Updating $(\{\beta_{n}\}_{n=1}^{N})^{(\tau+1)}$.
\end{itemize}

For \textbf{Step 1}, it is easy to find that the optimal solution of $\dot{\mathcal{T}}^{(\tau+1)}$ is
\begin{equation}
\label{equ16}
\dot{\mathcal{T}}^{(\tau+1)}=\mathcal{P}_{\Omega^{c}}
\left(\frac{1}{N}\sum_{n}^{N}\bigg({\rm{Fold}}_{n}(\dot{\mathbf{X}}_{[n]}^{(\tau)})-\frac{1}{\beta_{n}^{(\tau)}}{\rm{Fold}}_{n}(\dot{\mathbf{F}}_{[n]}^{(\tau)})\bigg)\right)+\dot{\mathcal{Y}},
\end{equation}
where $\Omega^{c}$ is the complement of $\Omega$, and we have used the fact that $\mathcal{P}_{\Omega^{c}}(\dot{\mathcal{Y}})=\mathbf{0}$ in (\ref{equ16}).

For \textbf{Step 2}, it can be decomposed into $N$ independent optimization problems, which can be solved in paralleled. For each $n$, $\dot{\mathbf{X}}_{[n]}^{(\tau+1)}$ is the optimal solution of the following
problem:
\begin{align}
\label{equ17}
\dot{\mathbf{X}}_{[n]}^{(\tau+1)}&=\mathop{{\rm{arg\, min}}}\limits_{\dot{\mathbf{X}}_{[n]}}\ \alpha_{n}\|\dot{\mathbf{X}}_{[n]}\|_{\ast}+\langle\dot{\mathbf{F}}_{[n]}^{(\tau)},\dot{\mathbf{T}}_{[n]}^{(\tau+1)}-\dot{\mathbf{X}}_{[n]}\rangle+\frac{\beta_{n}^{(\tau)}}{2}\|\dot{\mathbf{T}}_{[n]}^{(\tau+1)}-\dot{\mathbf{X}}_{[n]}\|_{F}^{2} \nonumber\\
&=\mathop{{\rm{arg\, min}}}\limits_{\dot{\mathbf{X}}_{[n]}}\ \alpha_{n}\|\dot{\mathbf{X}}_{[n]}\|_{\ast}+\frac{\beta_{n}^{(\tau)}}{2}\|\dot{\mathbf{X}}_{[n]}-(\dot{\mathbf{T}}_{[n]}^{(\tau+1)}+\frac{1}{\beta_{n}^{(\tau)}}\dot{\mathbf{F}}_{n}^{(\tau)})\|_{F}^{2} \nonumber\\
&=\mathop{{\rm{arg\, min}}}\limits_{\dot{\mathbf{X}}_{[n]}}\ \frac{\alpha_{n}}{\beta_{n}^{(\tau)}}\|\dot{\mathbf{X}}_{[n]}\|_{\ast}+\frac{1}{2}\|\dot{\mathbf{X}}_{[n]}-(\dot{\mathbf{T}}_{[n]}^{(\tau+1)}+\frac{1}{\beta_{n}^{(\tau)}}\dot{\mathbf{F}}_{n}^{(\tau)})\|_{F}^{2}.
\end{align}
Then, the optimal solution of (\ref{equ17}) can be obtained by the following theorem.
\begin{theorem}\footnote{Similar theorems for traditional real matrix case can be found in \cite{Furnival2017Denoising, DBLP:journals/tgrs/HeZZS16, DBLP:journals/tip/BenguaPTD17}, \emph{etc.}, we will demonstrate that the theorem still holds true for the quaternion matrix case in this paper.}
Let $\dot{\mathbf{Q}}\in\mathbb{H}^{N_{1}\times N_{2}}$ be a given quaternion matrix, then the QSVD (defined on \textbf{theorem 1}) of $\dot{\mathbf{Q}}$ with rank $r$ is
\begin{equation}
\label{ssvd}
\dot{\mathbf{Q}}=\dot{\mathbf{U}}\mathbf{\Lambda}\dot{\mathbf{V}}^{H}=\dot{\mathbf{U}}_{r}\mathbf{\Sigma}_{r}\dot{\mathbf{V}}_{r}^{H},
\end{equation}
where $\dot{\mathbf{U}}_{r}=[\dot{\mathbf{u}}_{1},\dot{\mathbf{u}}_{2},\ldots,\dot{\mathbf{u}}_{r}]\in\mathbb{H}^{N_{1}\times r}$ and $\dot{\mathbf{V}}_{r}=[\dot{\mathbf{v}}_{1},\dot{\mathbf{v}}_{2},\ldots,\dot{\mathbf{v}}_{r}]\in\mathbb{H}^{N_{2}\times r}$, $\mathbf{\Sigma}_{r}={\rm{diag}}(\sigma_{1},\ldots,\sigma_{r})$.
Define the quaternion matrix singular value thresholding operator
$\mathfrak{S}_{\xi}(\dot{\mathbf{Q}})=\dot{\mathbf{U}}_{r}\breve{\mathbf{\Sigma}}_{r}\dot{\mathbf{V}}_{r}^{H}$, where $\breve{\mathbf{\Sigma}}_{r}={\rm{diag}}\{{\rm{max}}(\sigma_{n}-\xi,0) \} (n=1,2,\ldots,r)$. Then the operator $\mathfrak{S}_{\xi}(\dot{\mathbf{Q}})$ obeys
\begin{equation}
\label{equ18}
\mathfrak{S}_{\xi}(\dot{\mathbf{Q}})=\mathop{{\rm{arg\, min}}}\limits_{\dot{\mathbf{X}}}\ \xi\|\dot{\mathbf{X}}\|_{\ast}+\frac{1}{2}\|\dot{\mathbf{X}}-\dot{\mathbf{Q}}\|_{F}^{2}.
\end{equation}
\end{theorem}
$\mathit{Proof.}$ It is obvious that the function $\mathfrak{N}(\dot{\mathbf{X}})=\xi\|\dot{\mathbf{X}}\|_{\ast}+\frac{1}{2}\|\dot{\mathbf{X}}-\dot{\mathbf{Q}}\|_{F}^{2}$ is strictly convex, hence there indeed exists a unique minimizer. We say that $\dot{\mathbf{X}}_{\ast}$ minimizes $\mathfrak{N}(\dot{\mathbf{X}})$ if and only if $\mathbf{0}$ is a subgradient of the function $\mathfrak{N}(\dot{\mathbf{X}})$ at the point $\dot{\mathbf{X}}_{\ast}$, \emph{i.e.},
\begin{equation}
\label{cftj}
\mathbf{0}\in\xi\nabla_{sub}\|\dot{\mathbf{X}}_{\ast}\|_{\ast}+\dot{\mathbf{X}}_{\ast}-\dot{\mathbf{Q}},
\end{equation}
where $\nabla_{sub}\|\dot{\mathbf{X}}_{\ast}\|_{\ast}$ denotes the subgradient set of $\|\dot{\mathbf{X}}_{\ast}\|_{\ast}$, which can be obtained by the following Lemma.
\begin{lemma}\cite{DBLP:journals/nla/JiaNS19} Suppose that $\dot{\mathbf{X}}\in\mathbb{H}^{N_{1}\times N_{2}}$ with rank $r$ has the QSVD as $\dot{\mathbf{A}}_{r}\mathbf{D}_{r}\dot{\mathbf{B}}_{r}^{H}=\sum_{n=1}^{r}d_{n}\dot{\mathbf{a}}_{n}\dot{\mathbf{b}}_{n}^{H}$, then
\begin{equation}
\label{subd}
\nabla_{sub}\|\dot{\mathbf{X}}\|_{\ast}=\{\sum_{n=1}^{r}\dot{\mathbf{a}}_{n}\dot{\mathbf{b}}_{n}^{H}+\dot{\mathbf{E}}\mid\dot{\mathbf{E}}\in\mathbb{H}^{N_{1}\times N_{2}},\  \dot{\mathbf{A}}_{r}^{H}\dot{\mathbf{E}}=\mathbf{0},\  \dot{\mathbf{E}}\dot{\mathbf{B}}_{r}=\mathbf{0}, \ \|\dot{\mathbf{E}}\|\leq 1\}.
\end{equation}	
\end{lemma}
We rewritten the QSVD of $\dot{\mathbf{Q}}$ (\emph{see} (\ref{ssvd})) as
\begin{equation*}
\dot{\mathbf{Q}}=\dot{\mathbf{U}}_{r}\mathbf{\Sigma}_{r}\dot{\mathbf{V}}_{r}^{H}=\sum_{n=1}^{r}\sigma_{n}\dot{\mathbf{u}}_{n}\dot{\mathbf{v}}_{n}^{H}=\sum_{n=1}^{r_{0}}\sigma_{n}\dot{\mathbf{u}}_{n}\dot{\mathbf{v}}_{n}^{H}+\sum_{n=r_{0}+1}^{r}\sigma_{n}\dot{\mathbf{u}}_{n}\dot{\mathbf{v}}_{n}^{H},
\end{equation*}
where $\sigma_{n}\leq \xi$ when $r_{0}+1 \leq n\leq r$. Then, setting $\dot{\mathbf{X}}_{\ast}=\mathfrak{S}_{\xi}(\dot{\mathbf{Q}})$, we have
\begin{equation*}
\dot{\mathbf{X}}_{\ast}=\sum_{n=1}^{r_{0}}(\sigma_{n}-\xi)\dot{\mathbf{u}}_{n}\dot{\mathbf{v}}_{n}^{H}+\sum_{n=r_{0}+1}^{r}0\dot{\mathbf{u}}_{n}\dot{\mathbf{v}}_{n}^{H},
\end{equation*}
and as a result
\begin{equation*}
\dot{\mathbf{X}}_{\ast}-\dot{\mathbf{Q}}=-\xi\left(\sum_{n=1}^{r_{0}}\dot{\mathbf{u}}_{n}\dot{\mathbf{v}}_{n}^{H}+\sum_{n=r_{0}+1}^{r}\frac{\sigma_{n}}{\xi}\dot{\mathbf{u}}_{n}\dot{\mathbf{v}}_{n}^{H}\right)
=-\xi\left(\sum_{n=1}^{r_{0}}\dot{\mathbf{u}}_{n}\dot{\mathbf{v}}_{n}^{H}+\dot{\mathbf{E}}\right),
\end{equation*}
where $\dot{\mathbf{E}}=\sum_{n=r_{0}+1}^{r}-\frac{\sigma_{n}}{\xi}\dot{\mathbf{u}}_{n}\dot{\mathbf{v}}_{n}^{H}$. According to \textbf{Lemma 1}, it is obvious that $\dot{\mathbf{X}}_{\ast}-\dot{\mathbf{Q}}\in-\xi\nabla_{sub}\|\dot{\mathbf{X}}\|_{\ast}$, \emph{i.e.}, when  $\dot{\mathbf{X}}_{\ast}=\mathfrak{S}_{\xi}(\dot{\mathbf{Q}})$, (\ref{cftj}) holds. Consequently, the $\mathfrak{S}_{\xi}(\dot{\mathbf{Q}})$ obeys the optimization problem (\ref{equ18}).

Therefore, based on \textbf{Theorem 2}, we can easily obtain the following optimization result of (\ref{equ17})
\begin{equation}
\label{equ19}
\dot{\mathbf{X}}_{[n]}^{(\tau+1)}=\mathfrak{S}_{\frac{\alpha_{n}}{\beta_{n}^{(\tau)}}}\left(\dot{\mathbf{T}}_{[n]}^{(\tau+1)}+\frac{1}{\beta_{n}^{(\tau)}}\dot{\mathbf{F}}_{n}^{(\tau)}\right).
\end{equation}

For \textbf{Step 3}, it can also be decomposed into $N$ independent optimization problems and be solved in paralled. For each $n$, $\dot{\mathbf{F}}_{[n]}^{(\tau+1)}$ is updated by the following equation:
\begin{equation}
\label{equ20}
\dot{\mathbf{F}}_{[n]}^{(\tau+1)}=\dot{\mathbf{F}}_{[n]}^{(\tau)}-\beta_{n}^{(\tau)}(\dot{\mathbf{X}}_{[n]}^{(\tau+1)}-\dot{\mathbf{T}}_{[n]}^{(\tau+1)}).
\end{equation}

For \textbf{Step 4}, since the dynamical $(\{\beta_{n}\}_{n=1}^{N})$
are usually preferred to speed up the convergence of the algorithm \cite{DBLP:journals/pami/HuZYLH13}, we use the following way to adaptive update $\beta_{n}$, for $n=1,2,\ldots,N$.
\begin{equation}
\label{equ21}
\beta_{n}^{(\tau+1)}={\rm{min}}(\beta_{n}^{max},\eta\beta_{n}^{(\tau)}),
\end{equation}
where $\beta_{n}^{max}$ is the default maximum of $\beta_{n}$, and $\eta\geq 1$ is a constant parameter. Generally, we set $\eta=\eta_{0}> 1$ if $\|\dot{\mathcal{T}}^{(\tau+1)}-\dot{\mathcal{T}}^{(\tau)}\|_{F}$ is small enough (\emph{e.g.}, 0.01), $\eta=1$ otherwise.

Finally, the proposed \textbf{L}ow-\textbf{R}ank \textbf{C}ompletion for \textbf{Q}uaternion \textbf{T}ensor (\textbf{LRC-QT}) method can be summarized in Table \ref{tab_algorithm}.
\begin{table}[htbp]
\caption{The low-rank completion for quaternion tensor (\textbf{LRC-QT}) method.}
\hrule
\label{tab_algorithm}
\begin{algorithmic}[1]
	\REQUIRE The observed quaternion tensor data $\dot{\mathcal{T}}\in\mathbb{H}^{N_{1}\times N_{2} \times\ldots \times N_{N}}$, the observed index set $\Omega$.
	\STATE \textbf{Initialize} $(\{\dot{\mathbf{X}}_{[n]}\}_{n=1}^{N})^{(0)}$, $(\{\dot{\mathbf{F}}_{[n]}\}_{n=1}^{N})^{(0)}$, $\{\alpha_{n}\}_{n=1}^{N}$, $(\{\beta_{n}\}_{n=1}^{N})^{(0)}$, $\{\beta_{n}^{max}\}_{n=1}^{N}$, $\eta_{0}$, $\epsilon$.
	\STATE \textbf{Repeat}
	\STATE Update $\dot{\mathcal{T}}^{(\tau+1)}$ using equation (\ref{equ16}).
	\STATE \% Lines 5 and 6 all can be performed in parallel.
	\STATE Update $(\{\dot{\mathbf{X}}_{[n]}\}_{n=1}^{N})^{(\tau+1)}$ using equation (\ref{equ19}).
	\STATE Update $(\{\dot{\mathbf{F}}_{[n]}\}_{n=1}^{N})^{(\tau+1)}$ using equation (\ref{equ20}).
	\STATE Update $\beta_{n}^{(\tau+1)}$ using equation (\ref{equ21}).
	\STATE  $\tau\longleftarrow \tau+1$.
	\STATE \textbf{Until} $\|\dot{\mathcal{T}}^{(\tau+1)}-\dot{\mathcal{T}}^{(\tau)}\|_{F}\leq\epsilon$.
	\ENSURE $\dot{\mathcal{T}}$.
\end{algorithmic}
\hrule
\end{table}

In the rest of this section, to facilitate direct processing of color image recovery issues, as a special case of the aforementioned method, we consider the low-rank quaternion matrix (second-order tensor) completion problem, \emph{i.e.},
\begin{equation}
\label{equ22}
\begin{split}
&\mathop{{\rm{minimize}}}\limits_{\dot{\mathbf{T}}}\quad \alpha\|\dot{\mathbf{T}}\|_{\ast}\\
&\text{subject to} \quad  \mathcal{P}_{\Omega}(\dot{\mathbf{T}})=\dot{\mathbf{Y}},
\end{split}
\end{equation}
where $\dot{\mathbf{Y}}\in\mathbb{H}^{M\times N}$ is the underlying complete quaternion matrix, $\dot{\mathbf{T}}\in\mathbb{H}^{M\times N}$ is the observed quaternion matrix. For color image recovery, model (\ref{equ22}) is  different from traditional matrix based and third-order tensor based models and is more advantageous than them. Since the traditional matrix based models are inherently developed for gray-level images. Although third-order tensor based models can deal with this problem, for this type of algorithms, the rank of a tensor is generally pretty hard to determine \cite{DBLP:journals/tip/ZhouLLZ18}, so they usually cannot offer the best low-rank approximation to a tensor. Besides, the tensor factorization or matricization based methods (\emph{see, e.g.}, \cite{DBLP:journals/pami/LiuMWY13}) are likely to destroy color pixel structure. In brief, the recovery theory for low-rank tensor completion problem is not
well established compared with that of matrix based completion problems \cite{DBLP:journals/nla/JiaNS19}.

For problem (\ref{equ22}), adding an additional variable quaternion matrix $\dot{\mathbf{X}}\in\mathbb{H}^{M\times N}$, we can obtain the following equivalent formulation:
\begin{equation}
\label{equ23}
\begin{split}
\mathop{{\rm{minimize}}}\limits_{\dot{\mathbf{T}}}\quad &\alpha\|\dot{\mathbf{X}}\|_{\ast}\\
\text{subject to} \quad &\dot{\mathbf{T}}=\dot{\mathbf{X}},\quad \mathcal{P}_{\Omega}(\dot{\mathbf{T}})=\dot{\mathbf{Y}}.
\end{split}
\end{equation}
Then we define the following augment Lagrangian function:
\begin{equation}
\label{equ24}
\mathcal{L}(\dot{\mathbf{T}}, \dot{\mathbf{X}}, \dot{\mathbf{F}},\beta)
=\alpha\|\dot{\mathbf{X}}\|_{\ast}+\langle\dot{\mathbf{F}},\dot{\mathbf{T}}-\dot{\mathbf{X}}\rangle+\frac{\beta}{2}\|\dot{\mathbf{T}}-\dot{\mathbf{X}}\|_{F}^{2},
\end{equation}
where $\beta$ is the penalty parameter, $\dot{\mathbf{F}}$ is the Lagrange multiplier. According to the ADMM framework, we independently update $\dot{\mathbf{T}}$, $\dot{\mathbf{X}}$, $\dot{\mathbf{F}}$, $\beta$, and we summarize the proposed \textbf{L}ow-\textbf{R}ank \textbf{C}ompletion for \textbf{Q}uaternion \textbf{M}atrix (\textbf{LRC-QM}) method in Table \ref{tab_algorithm2}.
\begin{table}[htbp]
\caption{The low-rank completion for quaternion matrix (\textbf{LRC-QM}) method.}
\hrule
\label{tab_algorithm2}
\begin{algorithmic}[1]
	\REQUIRE The observed quaternion matrix data $\dot{\mathbf{T}}\in\mathbb{H}^{M\times N}$, the observed index set $\Omega$.
	\STATE \textbf{Initialize} $\dot{\mathbf{X}}^{(0)}$, $\dot{\mathbf{F}}^{(0)}$, $\alpha$, $\beta^{(0)}$, $\beta^{max}$, $\eta_{0}$, $\epsilon$.
	\STATE \textbf{Repeat}
	\STATE $\dot{\mathbf{T}}^{(\tau+1)} \longleftarrow \mathcal{P}_{\Omega^{c}}
	\left(\dot{\mathbf{X}}^{(\tau)}-\frac{1}{\beta}\dot{\mathbf{F}}^{(\tau)}\right)+\dot{\mathbf{Y}}$ (where $\mathcal{P}_{\Omega^{c}}(\dot{\mathbf{Y}})=\mathbf{0}$).
	\STATE $\dot{\mathbf{X}}^{(\tau+1)} \longleftarrow\mathfrak{S}_{\frac{\alpha}{\beta^{(\tau)}}}\left(\dot{\mathbf{T}}^{(\tau+1)}+\frac{1}{\beta^{(\tau)}}\dot{\mathbf{F}}^{(\tau)}\right)$.
	\STATE $\dot{\mathbf{F}}^{(\tau+1)}\longleftarrow \dot{\mathbf{F}}^{(\tau)}-\beta^{(\tau)}(\dot{\mathbf{X}}^{(\tau+1)}-\dot{\mathbf{T}}^{(\tau+1)})$.
	\STATE $\beta^{(\tau+1)}\longleftarrow{\rm{min}}(\beta^{max},\eta\beta^{(\tau)})$ (we set $\eta=\eta_{0}> 1$ if $\|\dot{\mathbf{T}}^{(\tau+1)}-\dot{\mathbf{T}}^{(\tau)}\|_{F}$ is small enough (\emph{e.g.}, 0.01), $\eta=1$ otherwise).
	\STATE  $\tau\longleftarrow \tau+1$.
	\STATE \textbf{Until} $\|\dot{\mathbf{T}}^{(\tau+1)}-\dot{\mathbf{T}}^{(\tau)}\|_{F}\leq\epsilon$.
	\ENSURE $\dot{\mathbf{T}}$.
\end{algorithmic}
\hrule
\end{table}

\subsection{The computational complexity analysis}
For LRC-QT, the observed quaternion tensor data $\dot{\mathcal{T}}\in\mathbb{H}^{N_{1}\times N_{2} \times\ldots \times N_{N}}$. We assume, for simplicity, that $N_{1}= N_{2} =\ldots = N_{N}=I$. It is easy to see that the main per-iteration computational complexity lies in the update of $\dot{\mathbf{X}}_{[n]}, n=1,2,\ldots,N$, which requires computing $N$ QSVD of $I\times I^{N-1}$ quaternion matrices. There has been some quaternionic algorithms (\emph{e.g.}, \cite{Bihan2007Jacobi}) were proposed for computing the QSVD, nevertheless, they are too time-consuming. We propose to compute the QSVD using the isomorphism between $\mathbb{H}^{N_{1}\times N_{2}}$ and $\mathbb{C}^{2N_{1}\times 2N_{2}}$ (\emph{see} (\ref{definef})). According to (\ref{qsvd}), the computation of the QSVD for $N$ $I\times I^{N-1}$ quaternion matrices is equivalent to the computation of the SVD of $N$ $2I\times 2I^{N-1}$ complex matrices, which can be performed using well-established classical algorithms of SVD over $\mathbb{C}^{2N_{1}\times 2N_{2}}$ (the built-in function `svd' in MATLAB 2014b is used by us). We still assume that the computational complexity of SVD for a $K\times K$ complex matrix is about $\mathcal{O}(K^3)$.
Therefore, the whole computational complexity of LRC-QT for one iteration is about $\mathcal{O}(8NI^{2N-1})$. Analogously, For LRC-QM, the observed quaternion matrix data $\dot{\mathbf{T}}\in\mathbb{H}^{M\times N}$. Assuming $M=N=I$, the main per-iteration computational complexity lies in the update of $\dot{\mathbf{X}}$, which is about $\mathcal{O}(8I^{3})$.

\section{Simulation results}
\label{sec5}
In this section, simulations on some natural color videos and images are conducted to evaluate the performance of the proposed LRC-QT and LRC-QM methods. And we compare them with several existing state-of-the-art approaches, including  SiLRTC \cite{DBLP:journals/pami/LiuMWY13}, SPC \cite{DBLP:journals/tsp/YokotaZC16}, TMac (involving TMac-inc and TMac-dec) \cite{1930-8337_2015_2_601} and TCTF (which is not applicable in color video simulation) \cite{DBLP:journals/tip/ZhouLLZ18}.
All the simulations are run in MATLAB $2014b$ under Windows $7$ on a personal computer with $2.20$GHz CPU and $8.00$GB memory.

All color videos and images, in our simulations, are initially represented by fourth-order tensors $\mathcal{T}_{V}\in \mathbb{R}^{M\times N\times 3 \times T}$ and third-order tensors $\mathcal{T}_{I}\in \mathbb{R}^{M\times N\times 3}$ respectively. For LRC-QT, each color video is reshaped as a pure
quaternion tensor $\dot{\mathcal{T}}\in\mathbb{H}^{M\times N \times T}$ by using the following way:
\begin{equation*}
\dot{\mathcal{T}}=\mathcal{T}_{V}(:,:,1,:)i+\mathcal{T}_{V}(:,:,2,:)j+\mathcal{T}_{V}(:,:,3,:)k.
\end{equation*}
For LRC-QM, each color image is reshaped as a  pure
quaternion matrix $\dot{\mathbf{T}}\in\mathbb{H}^{M\times N}$ by using the following way:
\begin{equation*}
\dot{\mathbf{T}}= \mathcal{T}_{I}(:,:,1)i+\mathcal{T}_{I}(:,:,2)j+\mathcal{T}_{I}(:,:,3)k.
\end{equation*}
In addition, we uniformly generate the index set $\Omega$ at Gaussian random distribution, and define the sampling ratio (SR) as:
${\rm{SR}}=\frac{{\rm{numel}}(\Omega)}{M\times N\times 3 \times T}100\% \ (\text{or}\  {\rm{SR}}=\frac{{\rm{numel}}(\Omega)}{M\times N\times 3}100\%\ \text{for color images})$, where ${\rm{numel}}(\Omega)$  represents the number of observation entries in the index set $\Omega$.

\textbf{Quantitative assessment:} To evaluate
the performance of proposed methods, except visual quality, we employ three quantitative quality indexes, including the peak signal-to-noise ratio (PSNR), the structure similarity (SSIM) and the feature similarity (FSIM), which are respectively defined as: ${\rm{PSNR}}=10{\rm{log}}10\left( \frac{{\rm{Peakval}}^{2}}{{\rm{MSE}}}\right)$, where ${\rm{Peakval}}$ is taken from the range of the pixel value datatype (\emph{e.g.}, for uint8 pixel value, it is 255), ${\rm{MSE}}$ is the mean square error, \emph{i.e.}, ${\rm{MSE}}=\|\mathcal{X}-\mathcal{T}\|_{F}^{2}/{\rm{numel}}(\mathcal{X})$, where $\mathcal{X}$ and $\mathcal{T}$ are the recovered and truth data, respectively;
${\rm{SSIM}}=\frac{(2\mu_{\mathcal{T}}\mu_{\mathcal{X}}+C_{1})(2\sigma_{\mathcal{T}\mathcal{X}}+C_{2})}{(\mu_{\mathcal{T}}^{2}+\mu_{\mathcal{X}}^{2}+C_{1})(\sigma_{\mathcal{T}}^{2}+\sigma_{\mathcal{X}}^{2}+C_{2})}$, where $\mu_{\mathcal{T}}$, $\mu_{\mathcal{X}}$, $\sigma_{\mathcal{T}}$, $\sigma_{\mathcal{X}}$ and $\sigma_{\mathcal{T}\mathcal{X}}$ are the local means, standard deviations, and cross-covariance for images $\mathcal{T}$ and $\mathcal{X}$, $C_{1}=(0.01L)^{2}$,  $C_{2}=(0.03L)^{2}$,  $C_{3}=C_{2}/2$, $L$ is the specified dynamic range of the pixel values (average structure similarity index
over frames (ASSIM) is chosen for color videos);
${\rm{FSIM}}=\frac{\sum_{z\in \Delta}S_{L}(z)PC_{m}(z)}{\sum_{z\in \Delta}PC_{m}(z)}$,
where $\Delta$ demotes the whole image spatial domain. The phase
congruency for position $z$ of image $\mathcal{T}$  is denoted as $PC_{x}(\mathcal{T})$, then $PC_{m}(z)={\rm{max}}\{PC_{\mathcal{T}(z)},PC_{\mathcal{X}(z)}\}$, $S_{L}(z)$ is the gradient magnitude for position $z$ (average feature similarity index
over frames (AFSIM) is chosen for color videos).

For LRC-QT, we let $\mathcal{X}(:,:,1,:)={\rm{Imag}}_{1}(\dot{\mathcal{T}})$, $\mathcal{X}(:,:,2,:)={\rm{Imag}}_{2}(\dot{\mathcal{T}})$ and $\mathcal{X}(:,:,3,:)={\rm{Imag}}_{3}(\dot{\mathcal{T}})$, where ${\rm{Imag}}_{n}(\dot{\mathcal{T}})\:(n=1,2,3)$ denotes $n$-th imaginary part of $\dot{\mathcal{T}}$. Analogously, for LRC-QM, $\mathcal{X}(:,:,1)={\rm{Imag}}_{1}(\dot{\mathbf{T}})$, $\mathcal{X}(:,:,2)={\rm{Imag}}_{2}(\dot{\mathbf{T}})$ and $\mathcal{X}(:,:,3)={\rm{Imag}}_{3}(\dot{\mathbf{T}})$.

\textbf{Datasets:} In the simulations, we use the color video dataset: YUV Video Sequence\footnote{\url{http://trace.eas.asu.edu/yuv/}} where each sequence contains at least $150$ frames, and the color image dataset: Berkeley Segmentation Dataset (BSD)\footnote{\url{ https://www2.eecs.berkeley.edu/Research/Projects/CS/vision/bsds/}} which includes $300$ clean color images of size $481\times 321\times 3$.

\begin{figure*}[htbp]
\centering
\subfigure[]{\includegraphics[width=3.2cm,height=2cm]{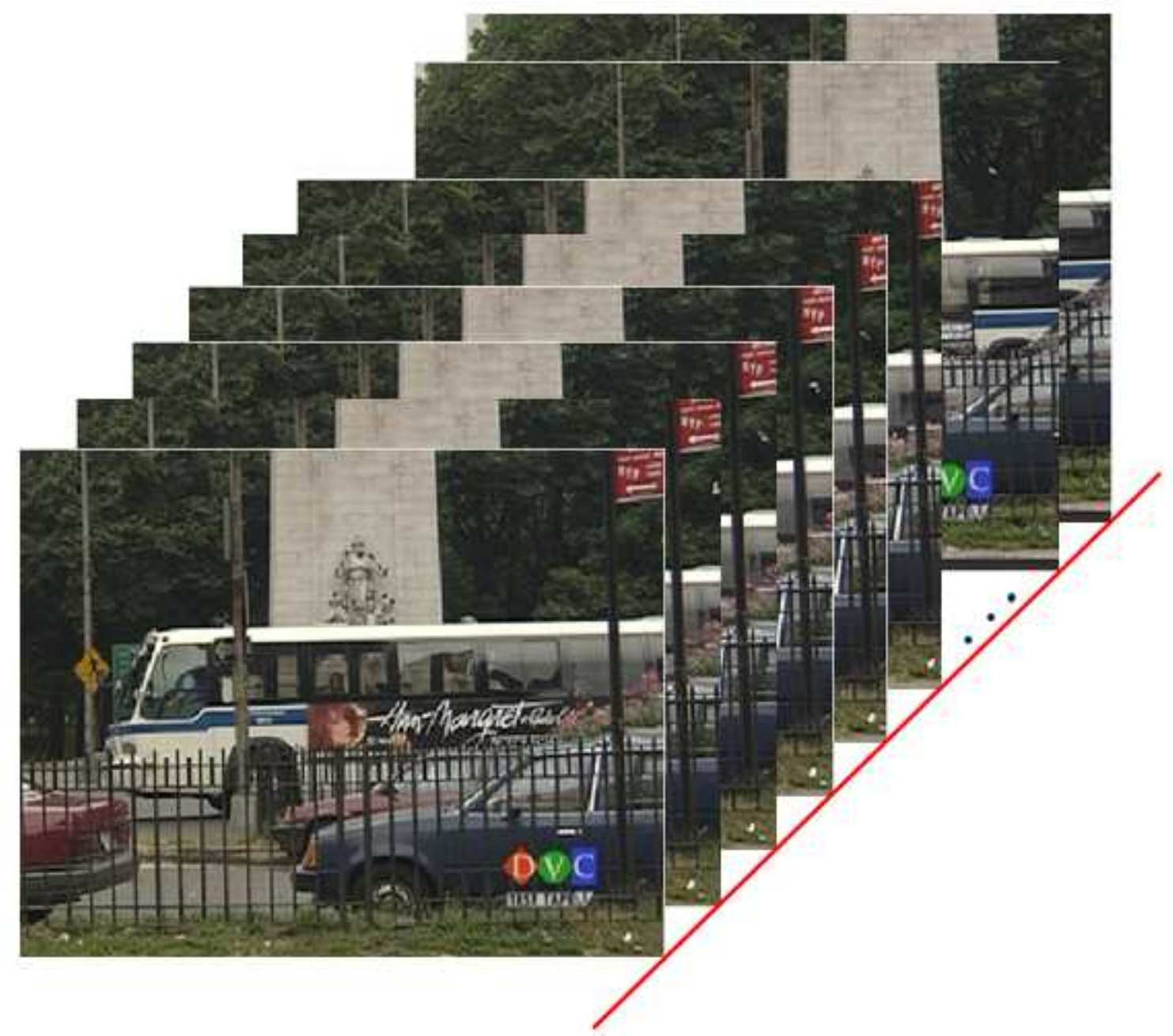}}\qquad
\subfigure[]{\includegraphics[width=2.8cm,height=2cm]{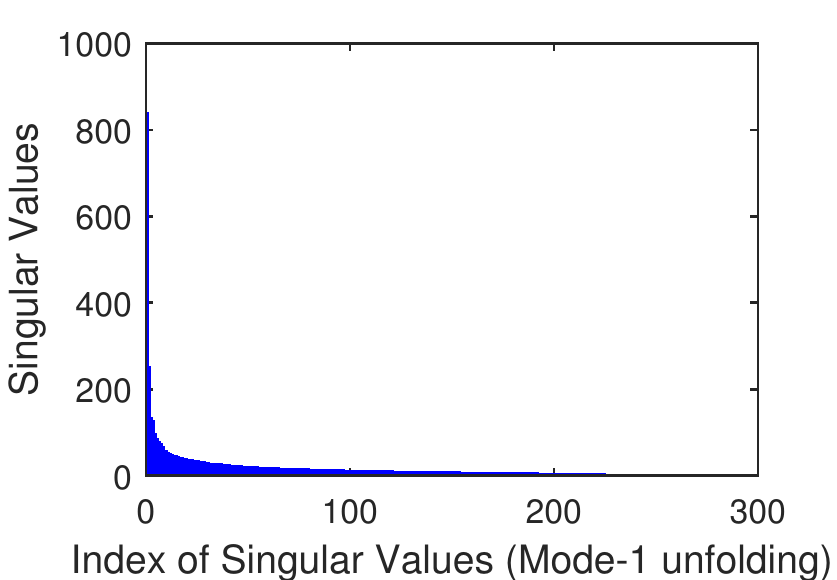}}
\subfigure[]{\includegraphics[width=2.8cm,height=2cm]{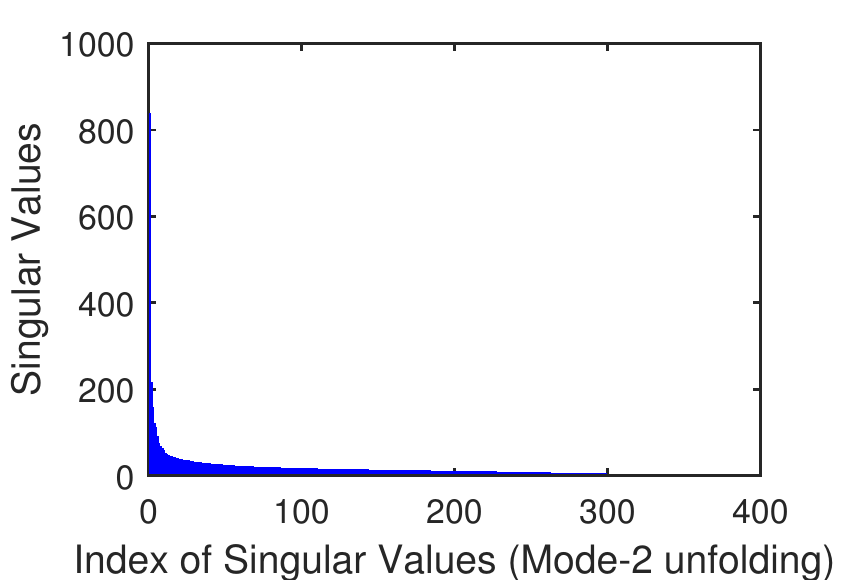}}
\subfigure[]{\includegraphics[width=2.8cm,height=2cm]{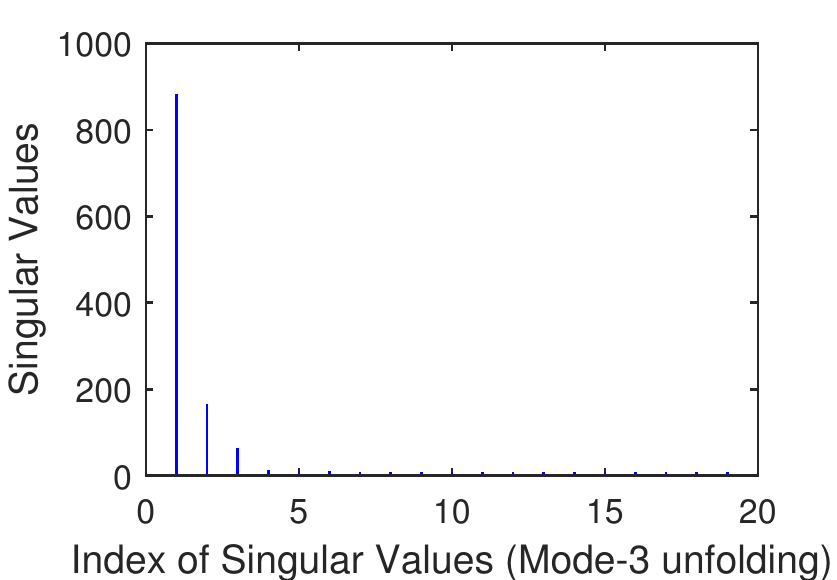}}\\	
\subfigure[]{\includegraphics[width=2.5cm,height=1.7cm]{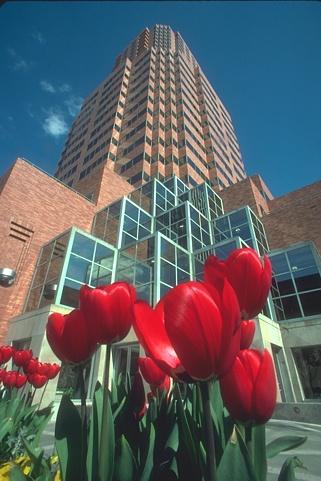}}\qquad
\subfigure[]{\includegraphics[width=2.8cm,height=2cm]{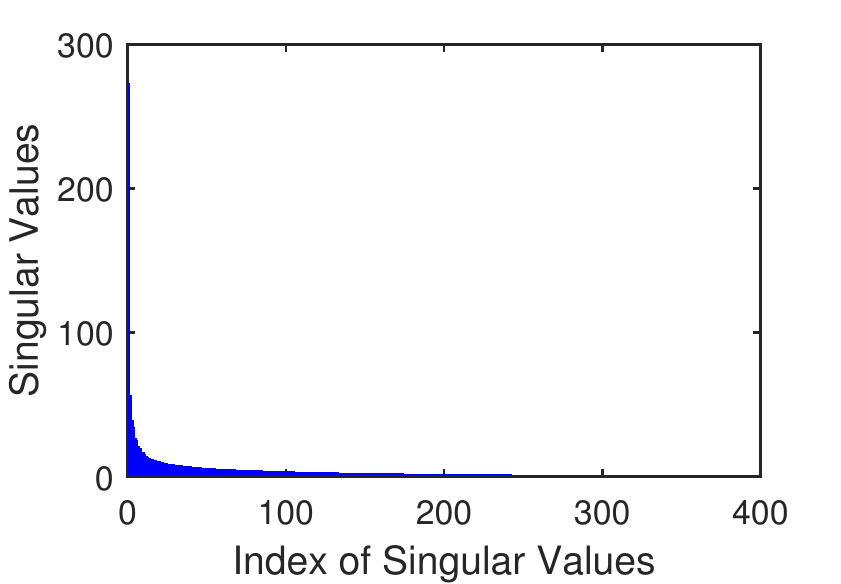}}
\caption{Illustration of the low-rank property of the color video (a) and color image (e). The (b), (c) and (d) respectively display the singular values of Mode-$1$, Mode-$2$ and Mode-$3$ unfolding quaternion matrices of (a). The (f) displays the singular values of (e).}
\label{fig_rank}
\end{figure*}

We first show that these color videos and color images can be well approximated and reconstructed by the low-rank quaternion tensors\footnote{The rank of quaternion tensor here refers to the Tucker rank defined by \textbf{Definition 5}.} and low-rank  quaternion matrices, respectively. As mentioned in \cite{4797640,DBLP:journals/tip/ZhouLLZ18}, \emph{etc.}, when the videos or images data are arranged into tensors or matrices, they lie on a union of low-rank subspaces approximately, which indicate the low-rank structure of the video or image data. This is also true for quaternion tensor and quaternion matrix data. For instance, in Figure.\ref{fig_rank} we display the singular values of one color video with size $288\times 352\times3\times20$ (reconstructed as third-order pure quaternion tensor with size $288\times 352\times20$) and one color image with size $481\times 321\times 3$ (reconstructed as pure quaternion matrix with size $481\times 321$), which are selected from the two datasets randomly. One can obviously see that most of the
singular values are very close to $0$ and much smaller than the first several larger singular values. So we could say that the color videos and color images can be well approximated by the low-rank quaternion tensors
and low-rank quaternion matrices, respectively, as we desired.

\textbf{Parameter and initialization settings:} For LRC-QT in Table \ref{tab_algorithm}, $(\{\dot{\mathbf{X}}_{[n]}\}_{n=1}^{3})^{(0)}=(\{\dot{\mathbf{F}}_{[n]}\}_{n=1}^{3})^{(0)}=(\mathbf{0}, \mathbf{0}, \mathbf{0})$, $\{\alpha_{n}\}_{n=1}^{3}=(2, 2, 10^{-3})$, $(\{\beta_{n}\}_{n=1}^{3})^{(0)}=(0.08, 0.08, 1)$, $\{\beta_{n}^{max}\}_{n=1}^{3}=(10^{3}, 10^{3}, 10^{3})$, $\eta_{0}=1.05$, $\epsilon=10^{-3}$. For LRC-QM in Table \ref{tab_algorithm2}, $\dot{\mathbf{X}}^{(0)}=\dot{\mathbf{F}}^{(0)}=\mathbf{0}$, $\alpha=2$, $\beta^{(0)}=0.08$, $\beta^{max}=10^{3}$, $\eta_{0}=1.05$, $\epsilon=10^{-3}$. For LRC-QT and LRC-QM, the parameters and initialization settings are just based on our experience and simulation results, and there may be better settings. For SiLRTC\footnote{\url{ http://www.cs.rochester.edu/ jliu/publications.html}}, SPC\footnote{\url{ http://ieeexplore.ieee.org/document/7502115/media}}, TMac\footnote{\url{http://www.caam.rice.edu/ yx9/TMac/}} and TCTF\footnote{\url{https://panzhous.github.io/}}, the codes of them are provided by their corresponding authors. The parameter settings and initialization methods of these algorithms are all based on the suggestions of their corresponding papers. Besides, for LRC-QT, as equations (\ref{equ19}) and (\ref{equ20}) show, we can update all $\dot{\mathbf{X}}_{[n]}$ and $\dot{\mathbf{F}}_{[n]}$ $(n=1,2,\ldots,N)$ parallelly, but for a fair comparison of the algorithm running time, we still employ the serial updating scheme in our code.

\begin{figure*}[htbp]
\centering
\subfigure[Bus]{\includegraphics[width=3.2cm,height=2.5cm]{rank1}}
\subfigure[News]{\includegraphics[width=3.2cm,height=2.5cm]{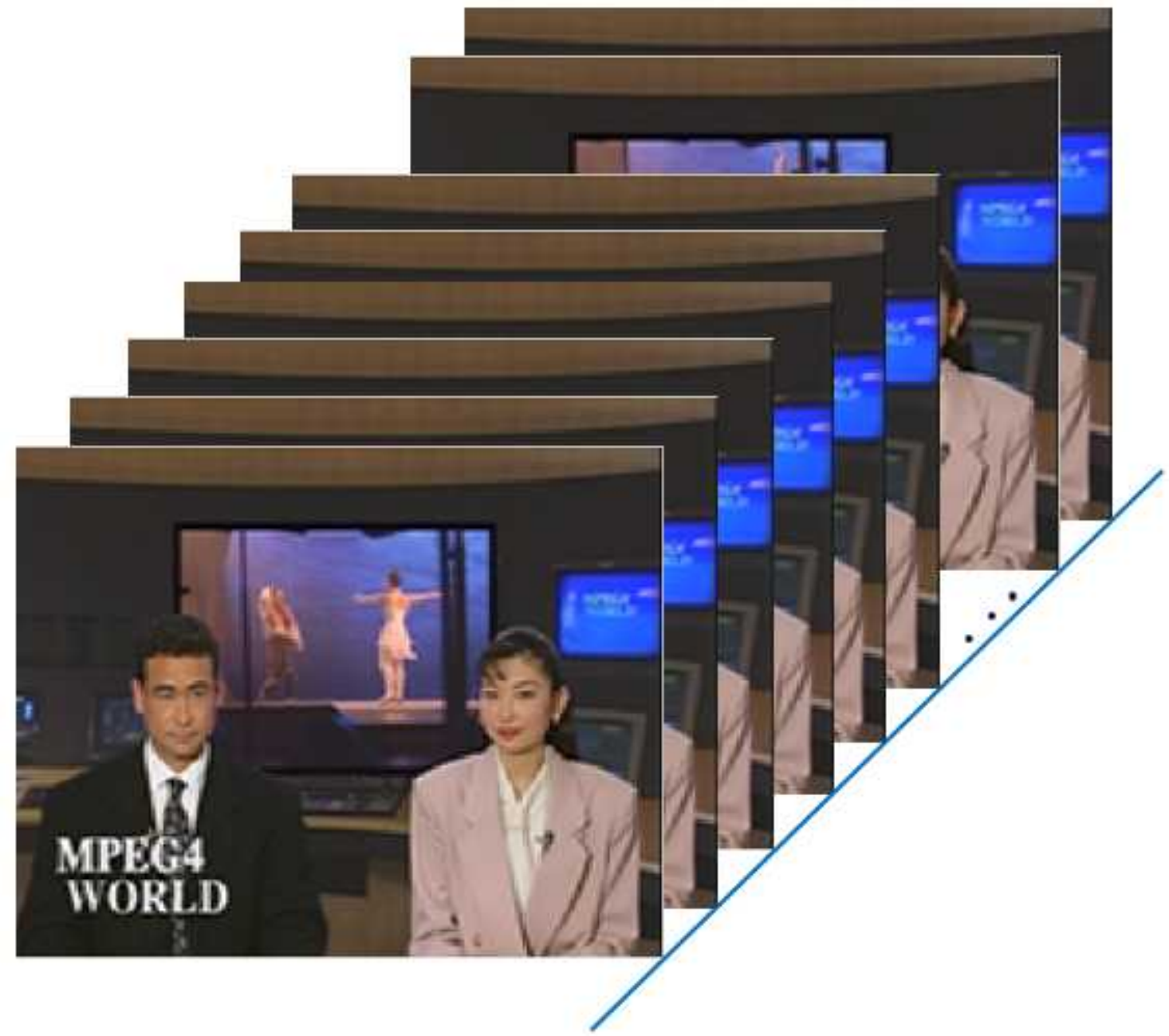}}
\subfigure[Salesman]{\includegraphics[width=3.2cm,height=2.5cm]{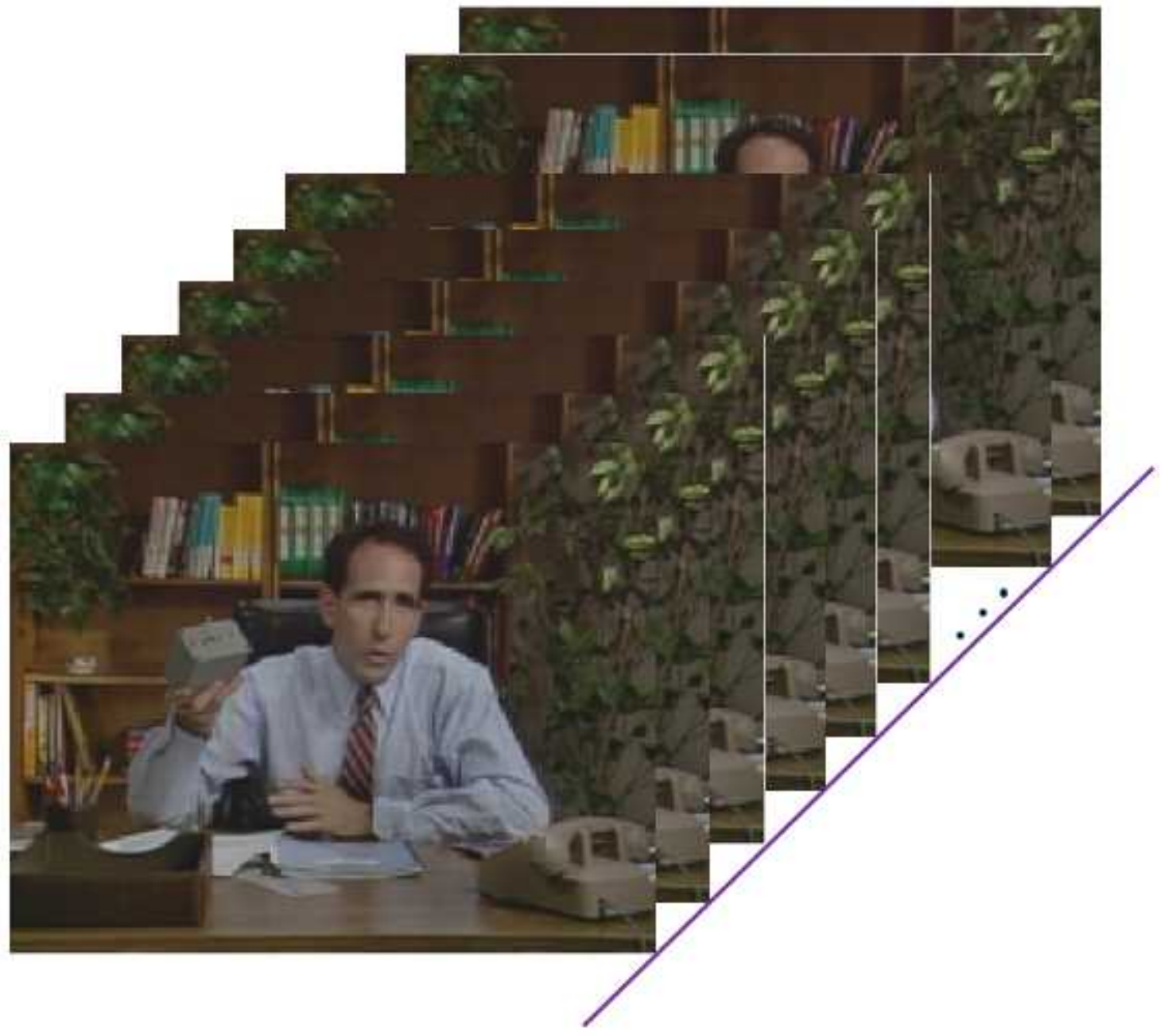}}
\subfigure[Suzie]{\includegraphics[width=3.2cm,height=2.5cm]{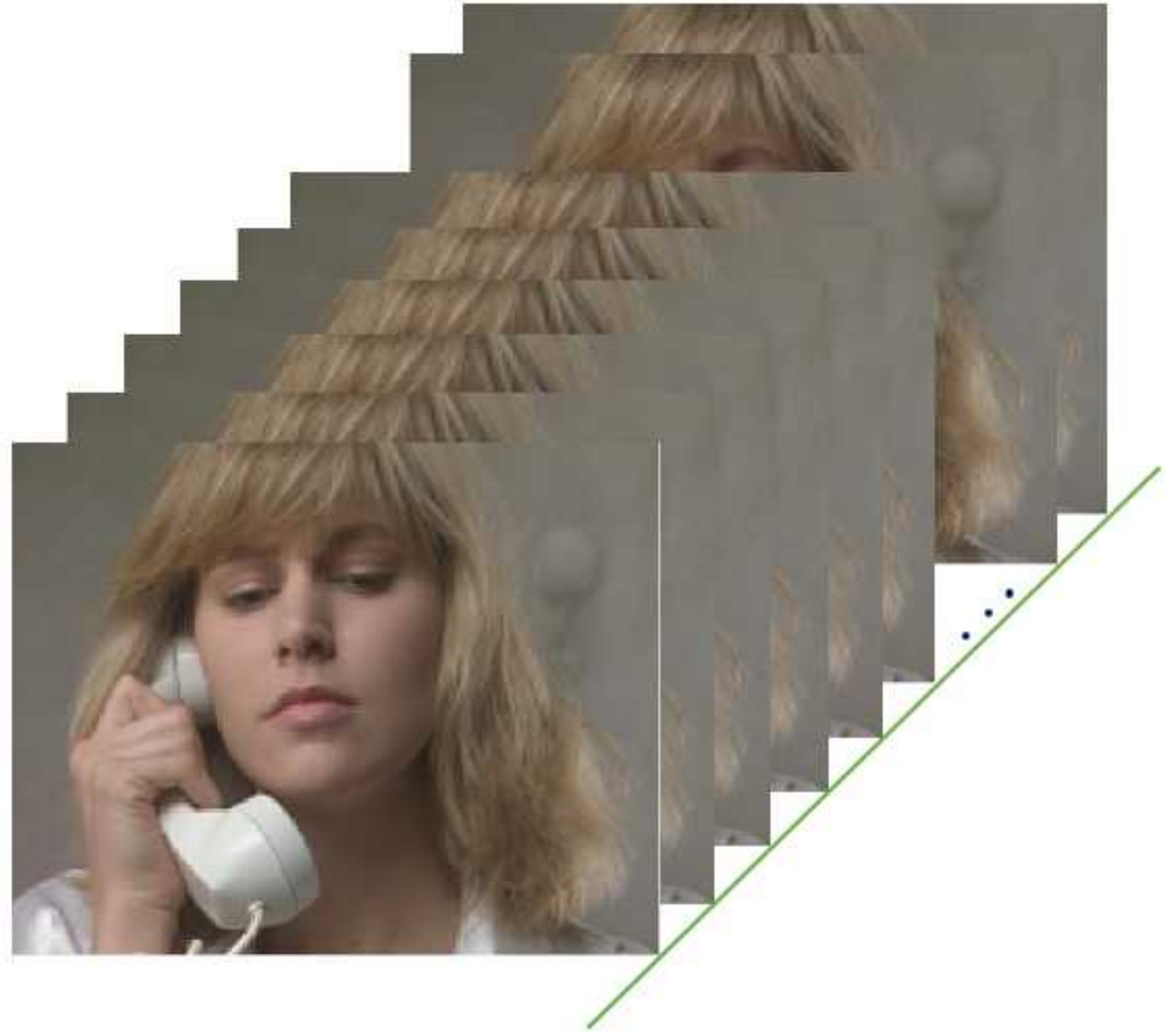}}
\caption{The four color videos (a) Bus, (b) News, (c) Salesman (d) Suzie, which are reconstructed by third-order pure quaternion tensors.}
\label{fig_zs}
\end{figure*}

\begin{figure*}[htbp]
\centering
\includegraphics[width=13cm,height=5.1cm]{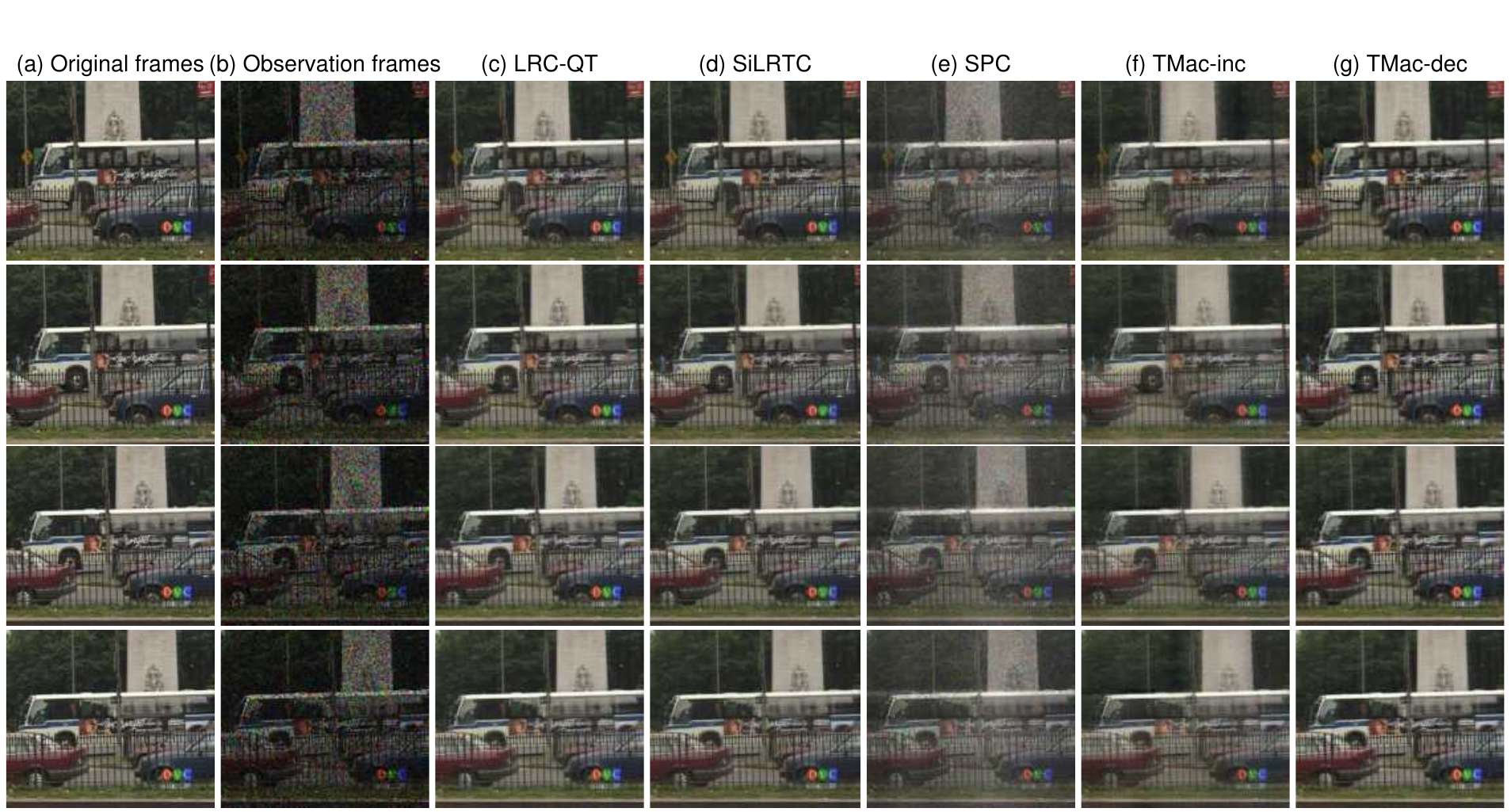}
\caption{The 1st, 8th, 15th and 20th frames (from top to
	bottom row) in the \textbf{Bus video}, with each column (from left to right) representing the original frames (a), observation frames with \textbf{${\rm{SR}}=\textbf{50\%}$} (b), recovery results of LRC-QT (c), SiLRTC (d),  SPC (e), TMac-inc (f) and TMac-dec (g).}
\label{busvideo}
\end{figure*}

\begin{figure*}[htbp]
\centering
\includegraphics[width=13cm,height=5.1cm]{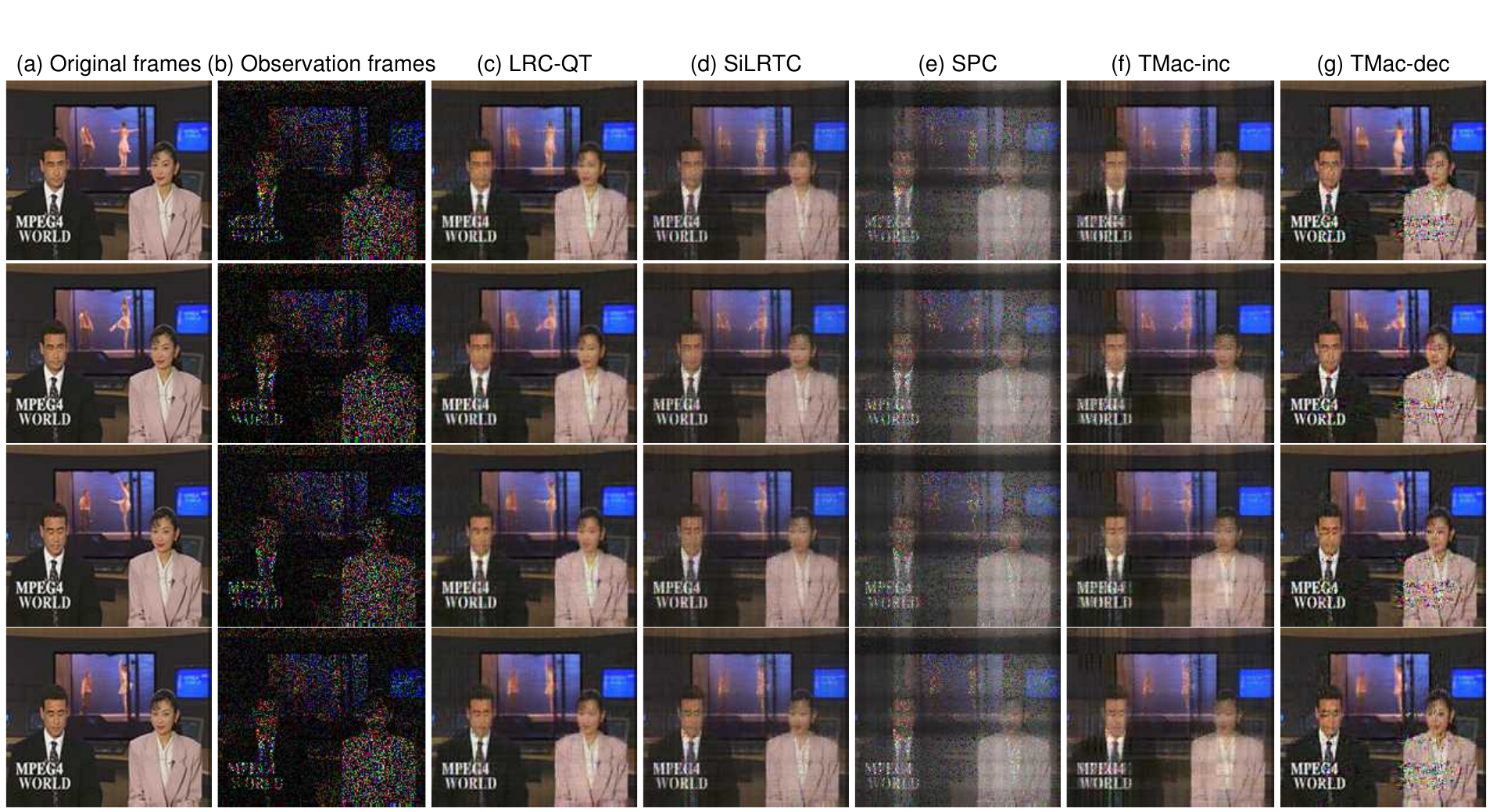}
\caption{The 1st, 8th, 15th and 20th frames (from top to
	bottom row) in the \textbf{News video}, with each column (from left to right) representing the original frames (a), observation frames with \textbf{${\rm{SR}}=\textbf{30\%}$} (b), recovery results of LRC-QT (c), SiLRTC (d),  SPC (e), TMac-inc (f) and TMac-dec (g).}
\label{newsvideo}
\end{figure*}

\begin{figure*}[htbp]
\centering
\includegraphics[width=13cm,height=5.1cm]{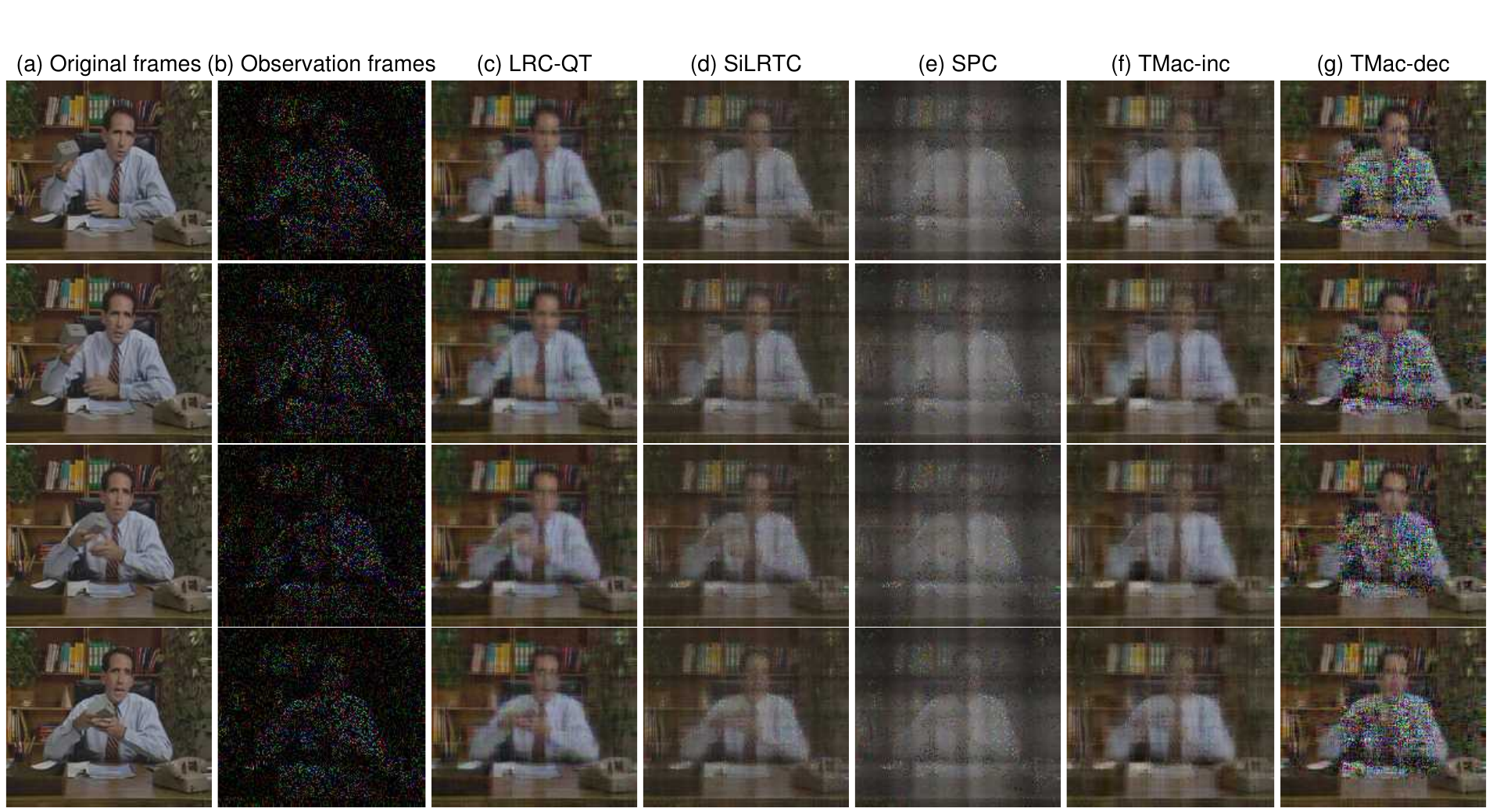}
\caption{The 1st, 8th, 15th and 20th frames (from top to
	bottom row) in the \textbf{Salesman video}, with each column (from left to right) representing the original frames (a), observation frames with \textbf{${\rm{SR}}=\textbf{20\%}$} (b), recovery results of LRC-QT (c), SiLRTC (d),  SPC (e), TMac-inc (f) and TMac-dec (g).}
\label{salesmanvideo}
\end{figure*}

\begin{figure*}[htbp]
\centering
\includegraphics[width=13cm,height=5.1cm]{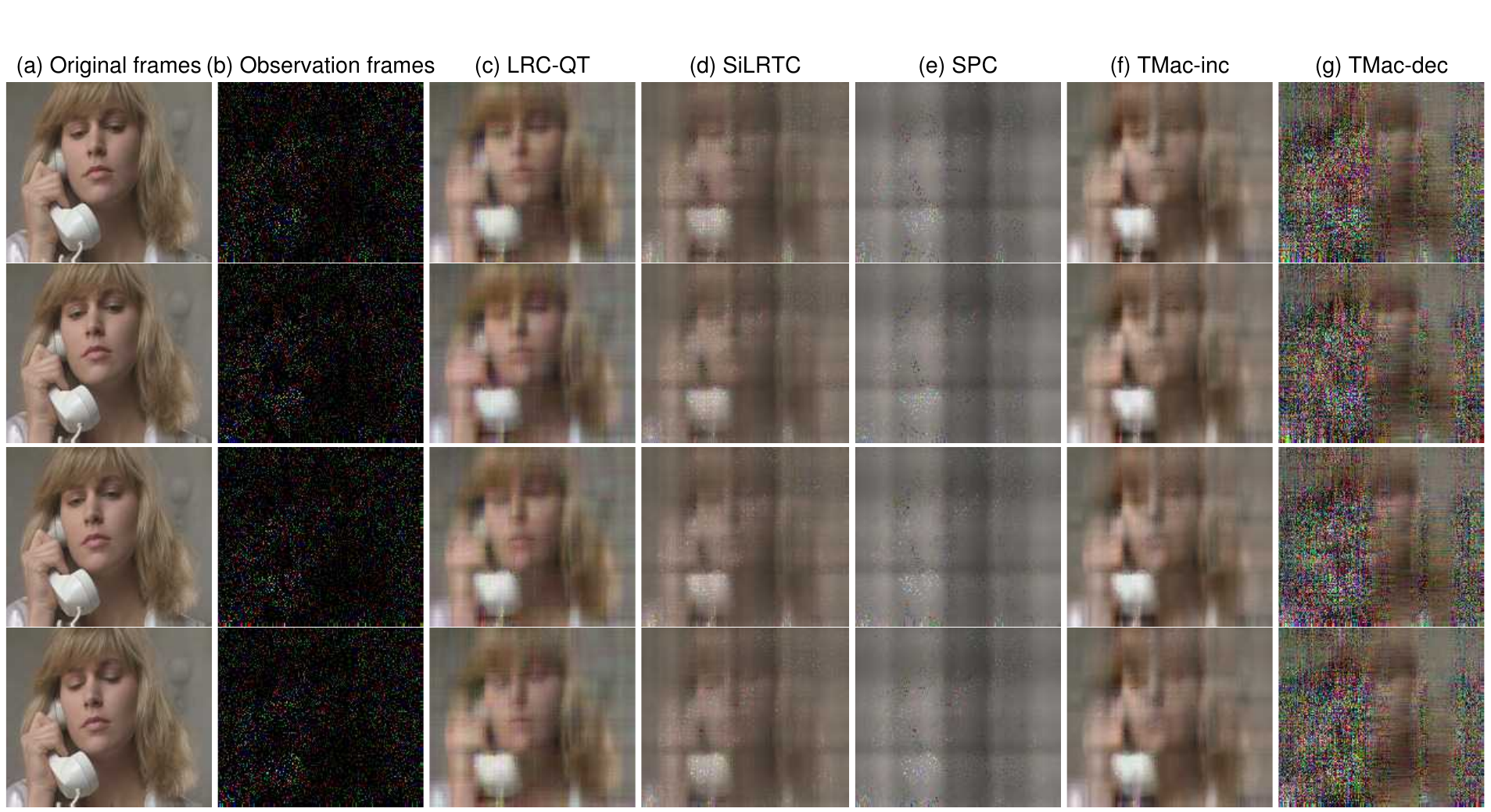}
\caption{The 1st, 8th, 15th and 20th frames (from top to
	bottom row) in the \textbf{Suzie video}, with each column (from left to right) representing the original frames (a), observation frames with ${\rm{SR}}=\textbf{10\%}$ (b), recovery results of LRC-QT (c), SiLRTC (d), SPC (e), TMac-inc (f) and TMac-dec (g).}
\label{suzievideo}
\end{figure*}

\textbf{Simulation 1:} In this simulation, we use four color videos (Bus, News, Salesman and Suzie) reconstructed by third-order pure quaternion tensors with size $288\times 253\times 20$, $144\times 176\times 20$, $144\times 176\times 20$ and $144\times 176\times 20$ respectively, and shown as Figure.\ref{fig_zs}, to evaluate the proposed LRC-QT for color video recovery. Figures.\ref{busvideo}-\ref{suzievideo} show the recovery results by different methods for the 1st, 8th, 15th and 20th frames of Bus video with   ${\rm{SR}}=50\%$, News video with  ${\rm{SR}}=30\%$, Salesman video with ${\rm{SR}}=20\%$ and Suzie video with ${\rm{SR}}=10\%$, respectively. We see from Figures.\ref{busvideo}-\ref{suzievideo} that the color video frames recovered by LRC-QT are visually better than those recovered by the other compared approaches. At very low SR (\emph{e.g.}, SR$=\%10$, \emph{see} Figure.\ref{suzievideo}), the advantage of LRC-QT seems to be more obvious. Table \ref{videotable} summaries the PSNR, ASSIM, AFSIM, iterations and running time of different methods on the four color videos with various SRs. From the results, one can find that for PSNR, ASSIM and AFSIM, LRC-QT reaches the highest values in most cases. We have reason to believe that this is mainly due to the advantage of quaternion representation of color pixel values. For the comparison of the number of iterations and running time, our approach is to compare them required by different compared methods when the PSNR obtained by all methods reaches a certain common pre-set value. We can see that although our method is not the fastest one, it can also achieve an acceptable PSNR value faster, and if we consider the fact that it can perform parallelly, the running time will be greatly reduced (Specifically, we can see that to achieve an acceptable PSNR value, our method requires only a few iterations, but the total time is not the shortest, which is closely related to the per-iteration computational complexity and how optimized the code is.). 
\begin{table}[p]
\caption{Quantitative quality indexes, iterations and running time (the display format is \underline{iterations/time}) of different methods on the four color videos with various SRs. The table is viewed better in zoomed PDF.}
\centering
\resizebox{15.5cm}{30mm}{
	\begin{tabular}{|c|c|ccccc|ccccc|ccccc|ccccc|c|}
		\hline
		\multicolumn{2}{|c|}{\diagbox{}{}} & &  &  \textbf{PSNR} & &  &  &  & \textbf{ASSIM} &  &  &  &  &\textbf{AFSIM}  &  &  & \multicolumn{6}{|c|}{\text{\textbf{Time (s)} (PSNR$\geq$ a certain value)}}\\
		\hline
		\textbf{Videos}& \textbf{SR} & LRC-QT& SiLRTC&SPC& TMac-inc& TMac-dec&LRC-QT& SiLRTC&SPC& TMac-inc& TMac-dec& LRC-QT& SiLRTC&SPC& TMac-inc& TMac-dec &LRC-QT& SiLRTC&SPC& TMac-inc& TMac-dec & PSNR$\geq \textbf{?}$ \\ \toprule
		\hline
		
		\multirow{5}{*}{Bus}
		& 50\% & \textbf{26.230} & 25.764 & 17.745 & 21.604 & 25.316 & \textbf{0.886} &0.877  & 0.561 & 0.725 & 0.841 & \textbf{0.997} & 0.996 & 0.978 & 0.988 & 0.992 &3/30.34  &37/55.55  & $\setminus$ &14/23.70  &8/39.92 &21 \\ \cline{2-23}
		
		& 40\% & \textbf{24.326} & 23.900 & 16.941 & 20.597 & 24.148 & \textbf{0.828} & 0.819 & 0.535 & 0.668 & 0.798 & \textbf{0.993} & 0.993 & 0.969 & 0.982 & 0.989 &4/40.58  &47/71.42  & $\setminus$ &16/26.64  &11/48.69 &20 \\ \cline{2-23}
		
		& 30\% &\textbf{22.696}  & 22.141 & 16.272 & 19.717 & 21.981 & \textbf{0.759} & 0.744 &0.477  & 0.607 &  0.745& \textbf{0.988} & 0.987 & 0.959 & 0.973 & 0.984 &5/43.94  &60/84.31  &$\setminus$ &18/27.99  &17/65.14 &19 \\ \cline{2-23}
		
		& 20\% &\textbf{21.805} &20.248  & 15.677 &18.876  & 21.367 &  \textbf{0.663} & 0.638 & 0.419 & 0.539 & 0.649 & \textbf{0.978} & 0.978 & 0.947 & 0.959 &0.977  &7/57.23  &86/117.48  & $\setminus$  &21/29.48  &37/124.01 & 18\\ \cline{2-23}
		
		& 10\% &\textbf{18.668}  & 17.863 & 15.145 & 18.074 & 17.304 & \textbf{0.527} & 0.482 & 0.387 & 0.463 & 0.471 & \textbf{0.959} & 0.958 & 0.930 & 0.941 & 0.959 & 9/72.63 &112/160.12  &22/18.29 &16/23.34  &65/207.62 & 15\\ \cline{2-23}
		\hline
		
		\multirow{5}{*}{News}
		& 50\% &32.252  &28.561  & 19.103 & 25.103 & \textbf{33.015} & \textbf{0.953} & 0.933 & 0.636 & 0.844 & 0.951 & \textbf{0.999} & 0.998 & 0.973 & 0.993 & 0.999 &6/9.77  &32/10.18  & $\setminus$ &20/8.72  &78/68.01 &25 \\ \cline{2-23}
		
		& 40\% &30.308  &26.308  &17.609  & 24.006 &\textbf{30.601}  & \textbf{0.941} & 0.903 & 0.584 & 0.809 &0.930  & \textbf{0.998} & 0.995 & 0.964 & 0.990 & 0.998 &7/11.45  &43/13.34  & $\setminus$ &22/9.43  &203/170.09 &24 \\ \cline{2-23}
		
		& 30\% & \textbf{27.153} & 24.578 & 16.630 & 23.001 & 22.024 & \textbf{0.893} & 0.861 & 0.538 & 0.774 & 0.813 & \textbf{0.996} & 0.993 & 0.953 & 0.986 & 0.992 &5/8.36  &54/16.13  & $\setminus$ &22/9.04  &982/798.06 &22 \\ \cline{2-23}
		
		& 20\% & \textbf{24.696} & 22.064 & 16.020 & 22.072 & 16.706 & \textbf{0.847} & 0.786 & 0.489 & 0.731 & 0.616 & \textbf{0.991} & 0.986 & 0.942 & 0.981 & 0.976 &4/6.68  &46/13.61  &17/4.05  &7/3.72 &988/805.69 &16 \\ \cline{2-23}
		
		& 10\% & \textbf{21.625} & 18.419 & 15,465 & 21.124 &15.711  & \textbf{0.756} & 0.628 & 0.335 & 0.674 & 0.532 & \textbf{0.980} & 0.970 & 0.924 & 0.974 &0.930  &7/10.97  &88/27.26  &17/4.69  &13/5.12  &991/831.61 &15 \\ \cline{2-23}
		\hline
		
		\multirow{5}{*}{Salesman}
		& 50\% &30.715 & 26.959 & 19.715 & 24.461 &\textbf{30.977}  &\textbf{0.959} &0.907  & 0.630 & 0.849 & 0.952 & \textbf{0.999} &0.996  & 0.982 & 0.990 & 0.998 &7/11.52  &27/10.95  &$\setminus$  &$\setminus$ &56/54.76 &25 \\ \cline{2-23}
		
		& 40\% &\textbf{29.311} & 25.205 & 18.906 & 23.363 & 29.113 & \textbf{0.937} &0.867  &0.588  &0.815  & 0.932 & \textbf{0.997} & 0.993 & 0.976 & 0.985 & 0.996 &7/11.57  &33/12.30  &$\setminus$  &15/6.48  &109/91.38 &23 \\ \cline{2-23}
		
		& 30\% & \textbf{26.324} &23.304  & 18.208 & 22.345 &25.011  &\textbf{0.884}  & 0.806 & 0.542 &0.771  & 0.862 & \textbf{0.993} & 0.989 & 0.967 & 0.977 & 0.978 &9/13.83  &46/16.21  &$\setminus$  &23/8.65  &524/492.81 &22 \\ \cline{2-23}
		
		& 20\% & \textbf{24.245} & 21.068 & 17.595 & 21.442 & 16.334 & \textbf{0.829} & 0.714 & 0.521 & 0.726 & 0.649 & \textbf{0.986} & 0.980 & 0.953 & 0.969 & 0.977 &4/6.88  &43/15.47  &9/2.38  &9/3.95  & $\setminus$ &17 \\ \cline{2-23}
		
		& 10\% & \textbf{21.224} & 17.797 & 17.018 &20.544  &12.204  &\textbf{0.717}  & 0.542 &0.443  & 0.672 & 0.400 & \textbf{0.969} & 0.959 & 0.929 & 0.958 & 0.921 &7/12.21 &106/32.48  &35/7.22  &21/6.94  &$\setminus$ &17 \\ \cline{2-23}
		\hline
		
		\multirow{5}{*}{Suzie}
		& 50\% & \textbf{34.801} & 31.302 & 21.749 & 28.953 & 33.048 &\textbf{0.980}  & 0.938 & 0.702 & 0.903 & 0.967 & \textbf{0.999} &  0.996& 0.967 & 0.991 & 0.998 &7/11.01  &34/11.43  & $\setminus$ &18/8.15  &280/255.55 &28 \\ \cline{2-23}
		
		& 40\% & \textbf{32.811} & 29.470 & 20.899 & 27.840 & 27.274 & \textbf{0.963} & 0.918 & 0.676 & 0.886 & 0.934 & \textbf{0.997} & 0.993 & 0.960 & 0.987 & 0.994 &9/13.89  &50/16.97  & $\setminus$ &20/8.23  &989/876.80 &27 \\ \cline{2-23}
		
		& 30\% & \textbf{30.613} & 27.484 & 20.272 & 26.810 & 22.179 & \textbf{0.944} & 0.886 & 0.633 & 0.872 & 0.871 & \textbf{0.993} & 0.989 & 0.952 & 0.981 & 0.988 &10/15.22  &61/19.37  & $\setminus$ &26/9.62 &$\setminus$ &26 \\ \cline{2-23}
		
		& 20\% & \textbf{28.247} & 24.989 & 19.519 & 25.894 & 17.629 & \textbf{0.915} & 0.834 & 0.541 & 0.851 & 0.721 & \textbf{0.987} & 0.980 & 0.945 & 0.975 & 0.973 &5/8.25  &60/17.55  &12/3.01  &10/3.16  &$\setminus$ &19 \\ \cline{2-23}
		
		& 10\% &\textbf{25.712} & 21.173  &18.881  & 25.033 & 15.103 &\textbf{0.837}  &0.738  & 0.490 & 0.834 & 0.482 &\textbf{0.974}  & 0.963 & 0.934 & 0.965 & 0.876 &9/15.05  &111/32.23  &21/4.64  &20/6.84  &$\setminus$ &18 \\ \cline{2-23}
		\hline
\end{tabular}}
\label{videotable}
\end{table}

\textbf{Simulation 2:} In this simulation, we use BSD dataset to evaluate the proposed LRC-QM for color image recovery. We randomly select 50 color images from this dataset reconstructed by pure quaternion matrices with size $481\times 321$.
\begin{figure*}[htbp]
\centering
\includegraphics[width=14cm,height=4cm]{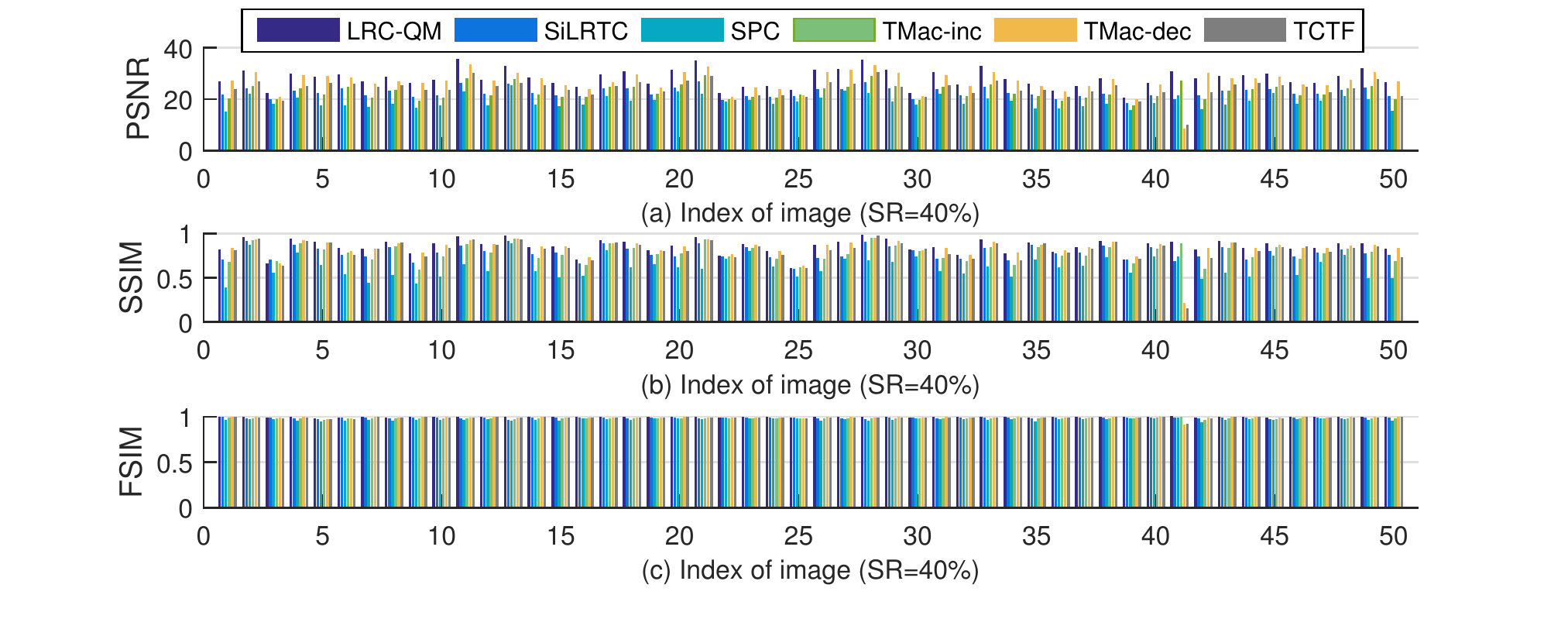}
\caption{Comparison of PSNR, SSIM and FSIM results of different algorithms for color image recovery on $50$ BSD images (${\rm{SR}}=40\%$). The figure is viewed better in zoomed PDF.}
\label{imag1}
\end{figure*}
\begin{figure*}[htbp]
\centering
\includegraphics[width=14cm,height=4cm]{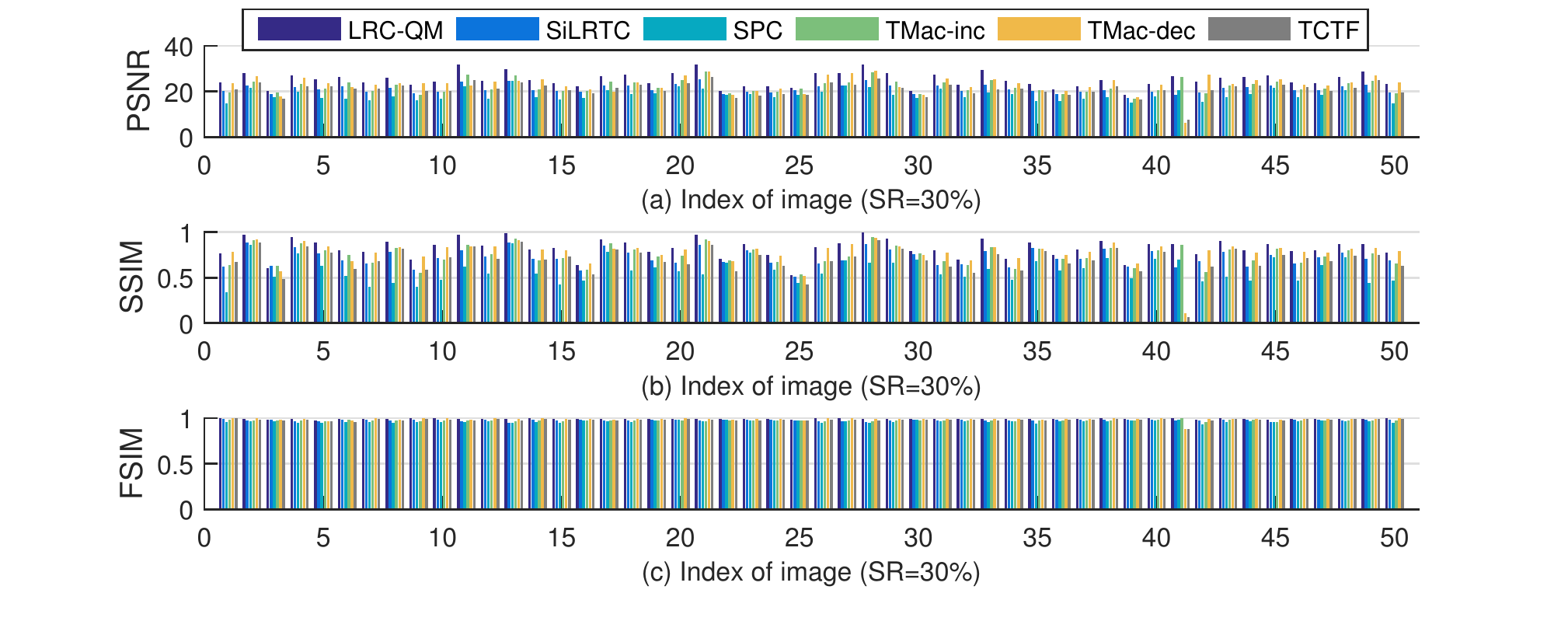}
\caption{Comparison of PSNR, SSIM and FSIM results of different algorithms for color image recovery on $50$ BSD images (${\rm{SR}}=30\%$). The figure is viewed better in zoomed PDF.}
\label{imag2}
\end{figure*}
\begin{figure*}[htbp]
\centering
\includegraphics[width=14cm,height=4cm]{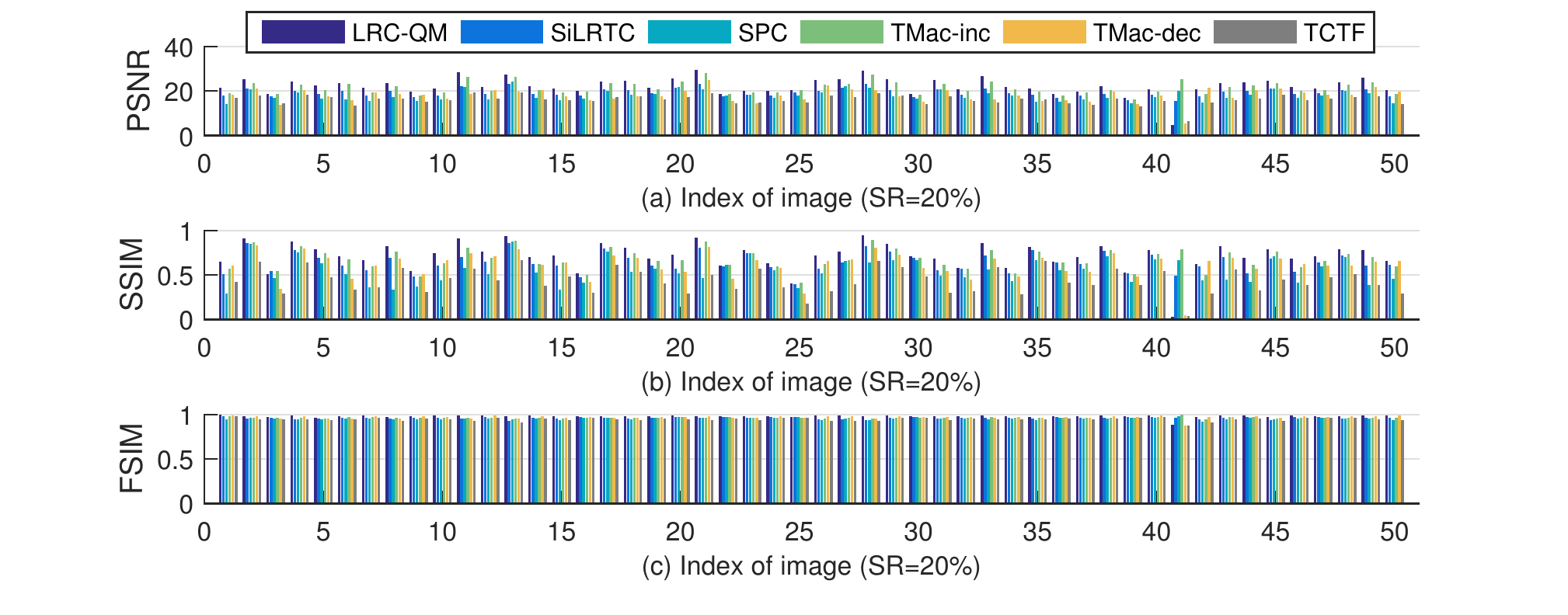}
\caption{Comparison of PSNR, SSIM and FSIM results of different algorithms for color image recovery on $50$ BSD images (${\rm{SR}}=20\%$). The figure is viewed better in zoomed PDF.}
\label{imag3}
\end{figure*}
In Figures.\ref{imag1}-\ref{imag3} we display the comparison of PSNR, SSIM and FSIM results of different methods for color image recovery on $50$ color images with various SRs. From the results, one can find that our LRC-QM approach performs better than the other methods in the vast majority of images. Besides the superiority on PSNR, the good performance of our method on SSIM and FSIM also demonstrates the advantages of the quaternion-based model.
\section{Conclusions}
\label{sec6}
Focusing on color videos and images recovery problems, this paper utilizing quaternions to represent the color pixels with RGB channels proposed a low-rank quaternion tensor completion method. Quaternion representation processes a color video or image holistically as a vector field and handles the coupling between the color channels naturally, and color information of source video or image is fully preserved. Although color video and images can also be represented as higher-order real tensors, the color structure will be destroyed during the process of matricization (\emph{e.g., mode-k unfolding}). Due to the special structure of the three imaginary parts of the quaternion, the relative positions of the three color channel pixels of one pixel are insensitive to the deformation of the quaternion tensor, which means that the color structure of the video can be completely maintained in the process of matricization. We adopted the ADMM framework to optimize the proposed model, which can guarantee the convergence of the algorithm. In addition, to facilitate direct processing of color image recovery issues, as a special case, we displayed the low-rank quaternion matrix completion model and optimization procedure separately. Theoretically, the proposed method can be well used to recover missing entries of any multidimensional data with color structures. In the simulation section, we mainly considered the color videos and images recovery problems. The results demonstrate the competitive performance (\emph{w.r.t.}, PSNR, SSIM and FSIM) of the proposed methods compared with several state-of-the-art approaches. 

Note that although the proposed method can well recover color videos and images, it needs to compute QSVD in each iteration, which is time-consuming and storage-intensive for large matrices. While the characteristic of less iteration of the algorithm can alleviate this shortcoming to some extent, in the future, we still aim to further explore better QSVD method to improve the efficiency of the algorithm, or to use some low-rank decomposition approaches to replace the nuclear norm minimization model.

\end{spacing}
\begin{spacing}{2}
\section*{Acknowledgments}
This work was supported by The Science and Technology Development Fund, Macau SAR (File no. FDCT/085/2018/A2).

\bibliographystyle{unsrt}
\bibliography{refs}

\begin{thebibliography}{10}

\bibitem{DBLP:journals/pr/LiWLT19}
Xiangrui Li, Andong Wang, Jianfeng Lu, and Zhenmin Tang.
\newblock Statistical performance of convex low-rank and sparse tensor
  recovery.
\newblock {\em Pattern Recognition}, 93:193--203, 2019.

\bibitem{DBLP:journals/pr/FanC17a}
Jicong Fan and Tommy W.~S. Chow.
\newblock Matrix completion by least-square, low-rank, and sparse
  self-representations.
\newblock {\em Pattern Recognition}, 71:290--305, 2017.

\bibitem{DBLP:journals/pr/FanC18}
Jicong Fan and Tommy W.~S. Chow.
\newblock Non-linear matrix completion.
\newblock {\em Pattern Recognition}, 77:378--394, 2018.

\bibitem{DBLP:journals/tip/ZouKW16}
Cuiming Zou, Kit~Ian Kou, and Yulong Wang.
\newblock Quaternion collaborative and sparse representation with application
  to color face recognition.
\newblock {\em {IEEE} Trans. Image Processing}, 25(7):3287--3302, 2016.

\bibitem{DBLP:journals/pr/QinJHLLZDZF19}
Binjie Qin, Mingxin Jin, Dongdong Hao, Yisong Lv, Qiegen Liu, Yueqi Zhu, Song
  Ding, Jun Zhao, and Baowei Fei.
\newblock Accurate vessel extraction via tensor completion of background layer
  in x-ray coronary angiograms.
\newblock {\em Pattern Recognition}, 87:38--54, 2019.

\bibitem{DBLP:journals/tip/ZhouLLZ18}
Pan Zhou, Canyi Lu, Zhouchen Lin, and Chao Zhang.
\newblock Tensor factorization for low-rank tensor completion.
\newblock {\em {IEEE} Trans. Image Processing}, 27(3):1152--1163, 2018.

\bibitem{DBLP:journals/sigpro/LongLCZ19}
Zhen Long, Yipeng Liu, Longxi Chen, and Ce~Zhu.
\newblock Low rank tensor completion for multiway visual data.
\newblock {\em Signal Processing}, 155:301--316, 2019.

\bibitem{6909886TNN}
Z.~{Zhang}, G.~{Ely}, S.~{Aeron}, N.~{Hao}, and M.~{Kilmer}.
\newblock Novel methods for multilinear data completion and de-noising based on
  tensor-svd.
\newblock In {\em 2014 IEEE Conference on Computer Vision and Pattern
  Recognition}, pages 3842--3849, June 2014.

\bibitem{DBLP:journals/siammax/KilmerBHH13}
Misha~Elena Kilmer, Karen~S. Braman, Ning Hao, and Randy~C. Hoover.
\newblock Third-order tensors as operators on matrices: {A} theoretical and
  computational framework with applications in imaging.
\newblock {\em {SIAM} J. Matrix Analysis Applications}, 34(1):148--172, 2013.

\bibitem{DBLP:journals/pami/ZhaoZC15}
Qibin Zhao, Liqing Zhang, and Andrzej Cichocki.
\newblock Bayesian {CP} factorization of incomplete tensors with automatic rank
  determination.
\newblock {\em {IEEE} Trans. Pattern Anal. Mach. Intell.}, 37(9):1751--1763,
  2015.

\bibitem{DBLP:journals/tsp/YokotaZC16}
Tatsuya Yokota, Qibin Zhao, and Andrzej Cichocki.
\newblock Smooth {PARAFAC} decomposition for tensor completion.
\newblock {\em {IEEE} Trans. Signal Processing}, 64(20):5423--5436, 2016.

\bibitem{DBLP:journals/pami/LiuMWY13}
Ji~Liu, Przemyslaw Musialski, Peter Wonka, and Jieping Ye.
\newblock Tensor completion for estimating missing values in visual data.
\newblock {\em {IEEE} Trans. Pattern Anal. Mach. Intell.}, 35(1):208--220,
  2013.

\bibitem{1930-8337_2015_2_601}
Yangyang Xu, Ruru Hao, Wotao Yin, and Zhixun Su.
\newblock Parallel matrix factorization for low-rank tensor completion.
\newblock {\em Inverse Problems and Imaging}, 9(2):601--624, 2015.

\bibitem{DBLP:journals/pr/YuWGGXP19}
Tingzhao Yu, Lingfeng Wang, Chaoxu Guo, Huxiang Gu, Shiming Xiang, and Chunhong
  Pan.
\newblock Pseudo low rank video representation.
\newblock {\em Pattern Recognition}, 85:50--59, 2019.

\bibitem{DBLP:journals/pr/ShiZWGZY19}
Jun Shi, Xiao Zheng, Jinjie Wu, Bangming Gong, Qi~Zhang, and Shihui Ying.
\newblock Quaternion grassmann average network for learning representation of
  histopathological image.
\newblock {\em Pattern Recognition}, 89:67--76, 2019.

\bibitem{DBLP:journals/mssp/GaiYW015}
Shan Gai, Guowei Yang, Minghua Wan, and Lei Wang.
\newblock Denoising color images by reduced quaternion matrix singular value
  decomposition.
\newblock {\em Multidim. Syst. Sign. Process.}, 26(1):307--320, 2015.

\bibitem{DBLP:journals/pr/HosnyD19}
Khalid~M. Hosny and Mohamed~M. Darwish.
\newblock New set of multi-channel orthogonal moments for color image
  representation and recognition.
\newblock {\em Pattern Recognition}, 88:153--173, 2019.

\bibitem{DBLP:journals/pr/LiLHS15}
Heng Li, Zhiwen Liu, Yali Huang, and Yonggang Shi.
\newblock Quaternion generic fourier descriptor for color object recognition.
\newblock {\em Pattern Recognition}, 48(12):3895--3903, 2015.

\bibitem{DBLP:journals/pr/ShaoSWCC14}
Zhuhong Shao, Huazhong Shu, Jiasong Wu, Beijing Chen, and Jean{-}Louis
  Coatrieux.
\newblock Quaternion bessel-fourier moments and their invariant descriptors for
  object reconstruction and recognition.
\newblock {\em Pattern Recognition}, 47(2):603--611, 2014.

\bibitem{DBLP:journals/ftml/BoydPCPE11}
Stephen~P. Boyd, Neal Parikh, Eric Chu, Borja Peleato, and Jonathan Eckstein.
\newblock Distributed optimization and statistical learning via the alternating
  direction method of multipliers.
\newblock {\em Foundations and Trends in Machine Learning}, 3(1):1--122, 2011.

\bibitem{DBLP:journals/tip/XuYXZN15}
Yi~Xu, Licheng Yu, Hongteng Xu, Hao Zhang, and Truong Nguyen.
\newblock Vector sparse representation of color image using quaternion matrix
  analysis.
\newblock {\em {IEEE} Trans. Image Processing}, 24(4):1315--1329, 2015.

\bibitem{doi:10.1080/14786444408644923}
Sir William Rowan Hamilton LL.D. P.R.I.A. F.R.A.S. Hon. M. R.~Soc. Ed. and
  Dub.~Hon. or~Corr.~M.
\newblock Ii. on quaternions; or on a new system of imaginaries in algebra.
\newblock {\em The London, Edinburgh, and Dublin Philosophical Magazine and
  Journal of Science}, 25(163):10--13, 1844.

\bibitem{Girard2007Quaternions}
Patrick~R. Girard.
\newblock {\em Quaternions, Clifford Algebras and Relativistic Physics}.
\newblock 2007.

\bibitem{ZHANG199721}
Fuzhen Zhang.
\newblock Quaternions and matrices of quaternions.
\newblock {\em Linear Algebra and its Applications}, 251:21--57, 1997.

\bibitem{DBLP:journals/siamrev/KoldaB09}
Tamara~G. Kolda and Brett~W. Bader.
\newblock Tensor decompositions and applications.
\newblock {\em {SIAM} Review}, 51(3):455--500, 2009.

\bibitem{doi:10.1137/0614071}
S.~Leurgans, R.~Ross, and R.~Abel.
\newblock A decomposition for three-way arrays.
\newblock {\em SIAM Journal on Matrix Analysis and Applications},
  14(4):1064--1083, 1993.

\bibitem{DBLP:journals/tip/BenguaPTD17}
Johann~A. Bengua, Ho~N. Phien, Hoang~Duong Tuan, and Minh~N. Do.
\newblock Efficient tensor completion for color image and video recovery:
  Low-rank tensor train.
\newblock {\em {IEEE} Trans. Image Processing}, 26(5):2466--2479, 2017.

\bibitem{DBLP:journals/ijon/TanCWZR14}
Huachun Tan, Bin Cheng, Wuhong Wang, Yu{-}Jin Zhang, and Bin Ran.
\newblock Tensor completion via a multi-linear low-n-rank factorization model.
\newblock {\em Neurocomputing}, 133:161--169, 2014.

\bibitem{DBLP:journals/siammax/GillisG11}
Nicolas Gillis and Fran{\c{c}}ois Glineur.
\newblock Low-rank matrix approximation with weights or missing data is
  np-hard.
\newblock {\em {SIAM} J. Matrix Analysis Applications}, 32(4):1149--1165, 2011.

\bibitem{Furnival2017Denoising}
Tom Furnival, Rowan~K. Leary, and Paul~A. Midgley.
\newblock Denoising time-resolved microscopy image sequences with singular
  value thresholding.
\newblock {\em Ultramicroscopy}, 178:112--124, 2017.

\bibitem{DBLP:journals/tgrs/HeZZS16}
Wei He, Hongyan Zhang, Liangpei Zhang, and Huanfeng Shen.
\newblock Total-variation-regularized low-rank matrix factorization for
  hyperspectral image restoration.
\newblock {\em {IEEE} Trans. Geoscience and Remote Sensing}, 54(1):178--188,
  2016.

\bibitem{DBLP:journals/nla/JiaNS19}
Zhigang Jia, Michael~K. Ng, and Guang{-}Jing Song.
\newblock Robust quaternion matrix completion with applications to image
  inpainting.
\newblock {\em Numerical Lin. Alg. with Applic.}, 26(4), 2019.

\bibitem{DBLP:journals/pami/HuZYLH13}
Yao Hu, Debing Zhang, Jieping Ye, Xuelong Li, and Xiaofei He.
\newblock Fast and accurate matrix completion via truncated nuclear norm
  regularization.
\newblock {\em {IEEE} Trans. Pattern Anal. Mach. Intell.}, 35(9):2117--2130,
  2013.

\bibitem{Bihan2007Jacobi}
Nicolas~Le Bihan.
\newblock Jacobi method for quaternion matrix singular value decomposition.
\newblock {\em Applied Mathematics and Computation}, 187(2):1265--1271, 2007.

\bibitem{4797640}
E.~J. {Candes} and B.~{Recht}.
\newblock Exact low-rank matrix completion via convex optimization.
\newblock In {\em 2008 46th Annual Allerton Conference on Communication,
  Control, and Computing}, pages 806--812, Sep. 2008.

\end{thebibliography}
\end{spacing}
\end{document}